\documentclass[a4paper,12pt]{amsart}
\usepackage{hyperref}
\usepackage{times}
\newtheorem{theorem}{Theorem}[section]
\newtheorem{lemma}[theorem]{Lemma}
\newtheorem{proposition}[theorem]{Proposition}
\newtheorem{corollary}[theorem]{Corollary}
\theoremstyle{definition}

\theoremstyle{remark}
\newtheorem{remark}[theorem]{Remark}
\numberwithin{equation}{section}
\DeclareMathOperator{\Div}{div}

\DeclareMathOperator{\curl}{curl}
\DeclareMathOperator{\Tr}{Trace}
\DeclareMathOperator{\supess}{sup\,ess}
\newcommand{\erre}{\mathbf{R}}
\newcommand{\D}{\mathcal{D}}
\newcommand{\F}{\mathcal{F}}
\newcommand{\E}{\mathbf{E}}
\newcommand{\Pb}{\mathbf{P}}
\newcommand{\e}{\mathtt{e}}
\newcommand{\N}{\mathbf{N}}
\newcommand{\auno}{$(A_1)$}
\newcommand{\adue}{$(A_2)$}
\newcommand{\atre}{$(A_3)$}
\newcommand{\aqua}{$(A_4)$}
\newcommand{\aquap}{$(A_4')$}
\newcommand{\acin}{$(A_5)$}
\newcommand{\asei}{$(A_6)$}
\newcommand{\sqn}{\sqrt{2\nu}}
\newcommand{\ox}{{\overline{x}}}
\newcommand{\ot}{{\overline{t}}}
\newcommand{\ou}{\overline{u}}
\newcommand{\ov}{\overline{v}}
\newcommand{\tv}{\tilde{v}}
\newcommand{\tD}{\tilde{\D}}
\newcommand{\tU}{\tilde{U}}
\newcommand{\s}{\sigma}
\newcommand{\ep}{\varepsilon}
\newcommand{\bs}{\mathcal B\kern-2.5pt{\scriptstyle\mathcal S}}
\newcommand{\ns}{\mathcal N\kern-5pt{\scriptstyle\mathcal S}}
\newcommand{\ualfa}{{\mathcal U}^\alpha(T)}
\newcommand{\tualfa}{{\mathcal U}^\alpha(\tau)}
\newcommand{\ualfam}{{\mathcal U}_M^\alpha(T)}
\newcommand{\tualfam}{{\mathcal U}_M^\alpha(\tau)}
\newcommand{\valfapi}{{\mathcal V}^{\alpha,p}(T)}
\newcommand{\tvalfapi}{{\mathcal V}^{\alpha,p}(\tau)}
\newcommand{\valfapil}{{\mathcal V}_L^{\alpha,p}(T)}
\newcommand{\tvalfapil}{{\mathcal V}_L^{\alpha,p}(\tau)}
\newcommand{\zero}[2]{\mathbf{0}_{#1\times #2}}
\newcommand{\set}[2]{\left\{\,#1\,\left|\,#2\,\right.\right\}}
\begin{document}
\hyphenation{di-par-ti-men-to ma-te-ma-ti-ca ap-pli-ca-ta
u-ni-ver-si-ta fi-ren-ze}
\title[A probabilistic representation for the vorticity...]{A
probabilistic representation for the vorticity of a 3D viscous
fluid and for general systems of parabolic equations}
\author[B. Busnello]{Barbara Busnello}
\address{Dipartimento di Matematica, Universit\`a di Pisa, via
         Buonarroti 2, 56127 Pisa, Italia}
\email{busnello@mail.dm.unipi.it}
\author[F. Flandoli]{Franco Flandoli}
\address{Dipartimento di Matematica Applicata, Universit\`a di 
         Pisa, Via Bonanno 25/b, 56126 Pisa, Italia}
\email{flandoli@dma.unipi.it}
\author[M. Romito]{Marco Romito}
\address{Dipartimento di Matematica, Universit\`a di Firenze, Viale 
         Morgagni 67/a, 50134 Firenze, Italia}
\email{romito@math.unifi.it}
\begin{abstract}
A probabilistic representation formula for general systems of linear
parabolic equations, coupled only through the zero-order term, is given. On
this basis, an implicit probabilistic representation for the vorticity
in a $3$D viscous fluid (described by the Navier-Stokes equations) is
carefully analysed, and a theorem of local existence and uniqueness is
proved.
\end{abstract}
\maketitle

\section{Introduction}

Consider the Navier-Stokes equation in $[0,T]\times\erre^3$
\begin{equation}\label{1.1}
\left\{\begin{array}{l}
\partial_t u+(u\cdot\nabla)u+\nabla p=\nu\Delta u+f \\
\Div u=0 \\
u(0,x)=u_0(x)
\end{array}\right.
\end{equation}
This equation describes, in \textsl{Eulerian coordinates}, the evolution
of a viscous incompressible Newtonian fluid, where $u$ is the velocity
field, $p$ the pressure, $f$ the body force and $\nu >0$ the kinematic
viscosity. The vorticity field $\xi =\curl u$ satisfies the equation
\begin{equation} \label{1.2}
\partial_t\xi+(u\cdot\nabla)\xi =\nu\Delta\xi+(\xi\cdot\nabla)u+g
\end{equation}
with $g=\curl f$. As we shall remark later on, the \textsl{stretching}
term $(\xi\cdot\nabla)u$ can be written in the form
$$
(\xi\cdot\nabla)u=\D_u\xi,
$$
where $\D_u=\frac12 (\nabla u+\nabla u^T)$, which better describes the
action of the deformation tensor $\D_u$ on $\xi$.  The analysis of the
vorticity field is a fundamental issue related to
questions like the possible emergence of singularities (see for instance
Beale, Kato and Majda \cite{BKM}, Constantin \cite{Con}), or the
description of  $3$D structures (see for instance Chorin \cite{Cho}).

The \textsl{Lagrangian} formulation of the fluid dynamics may be
important to analyse the vorticity field. Strictly speaking, the fluid
particles (we mean infinitesimal portions of fluid, not the single
molecules) move according to the deterministic law
$$
\dot X(t)=u(t,X(t)).
$$
However, a \textsl{virtual} Lagrangian dynamic of the particles of the form 
\begin{equation}\label{1.3}
dX(t)=u(t,X(t))\,dt+\sqn\,dW_t 
\end{equation}
(where $W_t$ is an auxiliary $3$D Brownian motion) allows us to describe
the evolution of quantities which are not only transported by the fluid,
but have a diffusive character. The vorticity has this property, as many
scalars or fields possibly  spreading into the fluid. Roughly speaking,
we prove the representation formula
$$
\xi(t,x)=\E[V(t,0)\xi_0(X(0))]+\int_0^t\E[V(t,s)g(s,X(s))]\,ds
$$
where $\E[\cdot]$ denotes the mean value with respect to the Wiener
measure, $\xi_0$ is the vorticity at time zero, $X(s)$ is the solution
of  equation \eqref{1.3} with final condition $X(t)=x$  and $V(r,s)$ is
the solution of the $3\times 3$ matrix equation
$$
\left\{\begin{array}{ll}
\frac{d}{dr}V(r,s)=\D_u(r,X(r))V(r,s),&\qquad r \in [s,t] \\
V(s,s)=I.
\end{array}\right.
$$

The present paper is devoted to explain the formula in detail, and use
it to prove a local-in-time existence and uniqueness result.  This paper
is in a sense the continuation of a paper of one of the authors (see
Busnello \cite{Bus}), where the $2$D case has been considered.  In the $2$D
case the stretching term $\D_u\xi$ is zero, so $V(r,s)=I$.  The
vorticity is purely transported and diffused, allowing for a
global-in-time control which yields global existence and uniqueness
results.  In Busnello \cite{Bus}, the probabilistic formula is used to
prove such a result, related to the deterministic work of Ben-Artzi
\cite{Ben}, following a suggestion of M. Friedlin.

Girsanov transformation is used in a basic part of the work, and the
Bismut-Elworthy formula is used to treat by probabilistic methods  also
the  Biot-Savart law, which reconstructs $u$ from $\xi$ (necessary to
solve \eqref{1.3}). In the $3$D case the Biot-Savart law and
its probabilistic representation are
$$
u(t,x)=-\frac1{4\pi}\int_{\erre^3}\frac{(x-y)\times\xi(y)}{|x-y|^3}\,dy
      =\frac12\int_0^\infty\frac1s\E[\xi(t,x+W_s)\times W_s]\,ds.
$$

In the present paper we extend as much as possible the probabilistic
approach of Busnello \cite{Bus} to the $3$D case.  Now, a
priori the stretching mechanism could produce singularities and blow-up
in $\xi$, so we can only work on a time interval $[0,\tau]$ depending on
the size of the data. This is the only possible result that also  the
analytic approaches to equation (\ref{1.1}) can reach at present.
Global existence for (\ref{1.1}) is known only at the level of weak
solutions, but we have to work at a higher level of regularity to deal
with the vorticity.  In certain function spaces, global existence (and
uniqueness) are known for sufficiently small data; in principle the
probabilistic formulation could  lead to such results, but we have found
some obstacles, so a probabilistic  proof of such a result remains an
open problem (except for the completely different approach of Le Jan and
Sznitman \cite{LeSz}).
\smallskip

The plan of the paper is the following. In Section \ref{snsrep} we state
the precise representation formula and the local existence and
uniqueness result for the Navier-Stokes equation, with the main lines of
its proof. However, the full proof of the representation formula and the
local result are based on three main items that we postpone to the next
three sections:
\begin{enumerate}
\item[\textit{(i)}] a general representation formula for linear systems
of parabolic equations, given in Section \ref{genpara};
\item[\textit{(ii)}] the probabilistic representation of Biot-Savart law and a
number of estimates on it, given in Section \ref{sbiotsavart};
\item[\textit{(iii)}] a series of estimates for the expected values
appearing in the formula for the vorticity, given in Section
\ref{lipschitz}.
\end{enumerate}
We have chosen this ordering to hi-light the results for the Navier-Stokes
equation at the beginning, for the reader who is not interested in the
long list of estimates and preliminaries necessary to prove the main
theorem. About item \textit{(i)} above, we remark that we use a method
due to Krylov (in the scalar case) that introduces new variables in
order to eliminate the zero order terms of the parabolic equation. Such
method in the case of systems coupled through the zero order part is
particularly interesting because it reduces the original system to a
decoupled one. The representation proved in Section \ref{genpara} can be
applied, in principle, to several other systems of equations appearing
in fluid dynamics, like the equation for $u$ itself (but the term
$\nabla p$ appears in the right-hand-side), the equation for the
magnetisation variable (see for example \cite{Cho}), the equation for
the transport of passive scalars.
\smallskip

Concerning the literature on the subject, at an advanced state of the
work we became aware of the interesting papers by Esposito, Marra,
Pulvirenti and Sciarretta \cite{EMPS} and by Esposito and Pulvirenti
\cite{EsPu} where a partially similar representation formula was
introduced; differently with respect to these papers we treat
probabilistically also the Biot-Savart law, we use different
probabilistic tools, we analyse in detail the general case of systems of
probabilistic equations to understand rigorously the equivalence with
the probabilistic representation and we prove the local existence and
uniqueness result in different function spaces (in particular, for a
class of less smooth initial conditions).

There is also a paper by Rapoport \cite{Rap} dealing with a general
class of equations on manifolds which in particular throw light on the
differential geometric structure of the formula.  Also the probabilistic
representation of systems of parabolic equations has been treated in the
literature under certain assumptions (our work seems to be more
general); see Kahane \cite{Kah} and Freidlin \cite{Fre}.

Finally, among the literature on probabilistic analysis of PDEs there
are possible connections with the geometric approach of Gliklikh
\cite{Gli}, with recent investigations on vortex method in $3$D by
Meleard and co-workers, by Giet \cite{Gie}, and more closely with a work
in preparation by Albeverio and Belopolskaya \cite{AlBe} where a
probabilistic representation for the velocity $u$ is
employed. Concerning the huge literature on the deterministic analysis
of the Navier-Stokes equations, results of local existence and
uniqueness have been proved in a great amount of function spaces, see
for instance collections of results in Cannone \cite{Can} and von Wahl
\cite{vWa}, or in many works of Kato, Solonnikov and many others. We
have not found a theorem exactly with the spaces used in the present
paper, but it may exist somewhere or it may be proved with an adaptation
of the existing techniques.

\subsection{A physical interpretation of the probabilistic formula for the 
vorticity}\label{interpreto}

\subsubsection{Evolution of the vorticity in the non-viscous case}

Let us first recall the well-known evolution of the vorticity field of
an incompressible \textit{non-viscous} fluid (therefore described by the
Euler equation). Let $\xi(t,x)$ be the value of the vorticity  at time
$t$ and point $x\in\erre^3$. The material point $x$ moves according to
the law
$$
\left\{\begin{array}{l}
\dot X(t)=u(t,X(t))\\
X(0)=x,
\end{array}\right.
$$
where $u$ is the velocity field of the fluid. From the Eulerian
description of the evolution of $\xi$
$$
\partial_t\xi+(u\cdot\nabla)\xi=\D_u\xi+g,
$$
we deduce the Lagrangian formulation 
\begin{equation}\label{9.2}
\frac{d}{dt}\xi(t,X(t))=\D_u(t,X(t))\xi(t,X(t))+g(t,X(t))
\end{equation}
which gives us
\begin{equation}\label{9.3}
\xi(t,X(t))=V(t,0)\xi_0(x)+\int _0^tV(t,s)g(s,X(s))\,ds
\end{equation}
where
\begin{equation}\label{9.4}
\left\{\begin{array}{ll}
\frac{d}{dr}V(r,s)=\D_u(r,X(r))V(r,s)&\qquad r\in[s,t]\\
V(s,s)=I.
\end{array}\right.
\end{equation}
Take $g=0$ for simplicity (the general case is similar); equations
(\ref{9.2}) and (\ref{9.3}) say that the initial vorticity $\xi_0(x)$ at
point $x$ is transported along the path $X(t)$, and during this motion it
is  modified by the deformation tensor. For instance, the vorticity is
stretched  when it is sufficiently aligned with the expanding directions
of $\D_u$; of course the relative position of $\xi$ with respect to the
expanding and  contracting (remember that $\Tr\D_u=0$) directions of $\D_u$
changes in time,  so $\xi(t,X(t))$ may undergo a complicate evolution
with stretching, rotations, contractions. Heuristic reasoning and
numerical experiments show a  predominance of the stretching mechanism,
and seem to indicate even a blow-up of $\xi(t,X(t))$ in finite time, for
certain initial point $x$.

If we want to know $\xi(\ot,\ox)$ at a certain time $\ot$ and point
$\ox$, we have to solve the backward equation
\begin{equation}\label{9.5}
\left\{\begin{array}{ll}
\dot X(t)=u(t,X(t))&\qquad t\in[0,\ot]\\
X(\ot)= \ox,
\end{array}\right.
\end{equation}
to find the initial position $x=X(0)$ which moves to $\ox$ at
time $\ot$; then 
\begin{equation} \label{9.6}
\xi(\ot,\ox)
= V^{\ot,\ox}(\ot,0)\xi_0(X^{\ot,\ox}(0))
 +\int_0^{\ot}V^{\ot,\ox}(\ot,s)g(s,X^{\ot,\ox}(s))\,ds
\end{equation}
where we have denoted by $X^{\ot,\ox}(\cdot)$  the solution of
(\ref{9.5}), to stress the dependence of the final condition $\ox$ at
time $\ot$, and by $V^{\ot,\ox}(r,s)$ the corresponding solution of
equation (\ref{9.4}).

\subsubsection{Path integral modification in the viscous case}

In the viscous case the position $X(t)$ of a material point still
evolves under the deterministic equation $\dot X(t)=u(t,X(t))$.  However,
the vorticity carried by the fluid particle at time $t=0$ is not simply
transported along its motion and modified by the action of the  tensor
$\D_u$; a diffusion of $\xi$ takes place.  Let us introduce a
\textsl{virtual} evolution of fluid particles, subject to $u$  and a
random diffusion:
\begin{equation}\label{9.7} 
dX(t)=u(t,X(t))\,dt+\sqn\,dW_t
\end{equation}
where $W_t$ is $3$D Brownian motion. Whether such a motion has a
physical meaning or not seems to be a similar question  to the case of
Feynman paths in Feynman integrals approach to quantum physics.  The
initial vorticity $\xi_0(x)$ decomposes, in a sense, in infinitesimal
components along the different random solutions of (\ref{9.7}),
proportional  to the probability of each evolution (strictly speaking
such probabilities  are zero).  If $X(t,\omega)$ is a path given by
(\ref{9.7}), let us denote by $P(\omega)$ its probability, ignoring for
a moment that $P(\omega)=0$; then an amount of vorticity equal to
$\xi_0(x)P(\omega)$ travels along $X(t,\omega)$ and is subject to the
action $V(t,s)$ of $\D_u$ along the path:
$$
\xi(t,X(t,\omega))P(\omega)=V(t,0)\xi_0(x)P(\omega)+\int _0^t
V(t,s) g(s,X(s, \omega))]\,ds\,P(\omega)
$$
(the reasoning for the integral effect of $g$ is similar, and we omit
it).  Now $\xi(t,X(t,\omega))P(\omega)$ is not the total value of the
vorticity field at time $t$ and point $\ox=X(t,\omega)$, but it is only
the contribution due to the $\omega$-evolution  started from position
$x$: other initial positions and other evolutions will reach the point
$\ox$ at time $t$, and we have to add all these contributions.
Therefore to compute $\xi(\ot,\ox)$ at a certain time  $\ot$ and point
$\ox$ we have to solve the backward stochastic equation
$$
\left\{\begin{array}{ll}
dX(t)=u(t,X(t))\,dt+\sqn\,dW_t&\qquad t\in[0,\ot] \\
X(\ot)=\ox
\end{array}\right.
$$
to find the various positions $X(0,\omega)$ which move to $\ox$ at  time
$\ot$ under different noise paths $W(t,\omega)$; at the  heuristic level
each $\omega$ gives a contribution $\xi(\ot,\ox; \omega))P(\omega)$ to
$\xi(\ot,\ox)$ given by
$$
\xi(\ot,\ox;\omega))P(\omega)
= V^{\ot,\ox}(\ot,0;\omega)\xi_0(X^{\ot,\ox}(0;\omega))P(\omega)
 +\int _0^{\ot}V^{\ot,\ox}(\ot,s)g(s,X^{\ot,\ox}(s;\omega))\,ds\,P(\omega)
$$
(see (\ref{9.6}) and (\ref{9.2})), so the total $\xi(\ot,\ox)$ is given by
$$
\xi(\ot,\ox)
= \E[V^{\ot,\ox}(\ot,0)\xi_0(X^{\ot,\ox}(0))]
 +\int _0^{\ot}\E[V^{\ot,\ox}(\ot,s)g(s,X^{\ot,\ox}(s))]\,ds.
$$
This is the heuristic derivation and the physical explanation of the
formula.


\section{Main result on the probabilistic representation for the
vorticity}\label{snsrep} 

\subsection{Some definitions and notations}

First we recall some classical spaces, like the space
$L^p(\erre^3,\erre^3)$ of 3D vector fields whose $p$-power is summable,
with norm
$$
\|f\|_p=\bigl(\int_{\erre^3}|f(x)|^p\,dx\bigr)^{\frac1p},
$$
the space $C^k_b(\erre^3,\erre^3)$ of k-times differentiable vector
fields, with norm
$$
\|g\|_{C^k_b}=\sum_{|\beta|\le k}\|D^\beta g\|_\infty
$$
and finally the space $C^{k,\alpha}_b(\erre^3,\erre^3)$ of vector fields
whose $k^\text{th}$-order derivatives are H\"older-continuous with
exponent $\alpha$, with norm
$$
\|g\|_{C^{k,\alpha}_b}=\|g\|_{C^k_b}+[g]_{k+\alpha},
$$
where
$$
[g]_{k+\alpha}=\sum_{|\beta|=k}\sup_{x,y\in\erre^3}\frac{|D^\beta
g(x)-D^\beta g(y)|}{\|x-y\|^\alpha}.
$$
Next we define the spaces where our problem will be set. The velocity
field of Navier-Stokes equations will be in the space
\begin{equation}\label{ualfa}
\ualfa=\set{u\in C([0,T];C^1_b(\erre^3,\erre^3))\cap
L^\infty(0,T;C^{1,\alpha}_b(\erre^3,\erre^3))}{\Div u(t)=0},
\end{equation}
endowed with the norm
$$
\supess_{0\le t\le T}\|u(t)\|_{C^{1,\alpha}_b},
$$
while the vorticity will be in the space
\begin{equation}\label{valfapi}
\valfapi=C\bigl([0,T];C_b(\erre^3,\erre^3)\cap
L^p(\erre^3,\erre^3)\bigr)\cap L^\infty\bigl(0,T;C^\alpha_b(\erre^3,\erre^3)\bigr),
\end{equation}
endowed with the norm
$$
\supess_{0\le t\le T}\|v(t)\|_{L^p\cap C^\alpha_b},
$$
where $\|\psi\|_{C^\alpha_b\cap L^p}=\|\psi\|_p+\|\psi\|_{C^\alpha_b}$.
We will use also the space
\begin{equation}\label{ualfam}
\ualfam=\set{u\in\ualfa}{\supess_{t\le T}\|u(t)\|_{C^{1,\alpha}_b}\le M},
\end{equation}
and the space
\begin{equation}\label{valfapil}
\valfapil=\set{\psi\in\valfapi}{\supess_{t\le T}\|\psi(t)\|_{L^p\cap C^\alpha_b}\le L}.
\end{equation}
\subsection{Probabilistic representation for the vorticity}
The formulation of the three dimensional Navier-Stokes equations
\begin{eqnarray}\label{nsformulation}
&&\partial_t u-\nu\Delta u+(u\cdot\nabla)u+\nabla P=f,\notag\\
&&\Div u=0,\\
&&u(0,x)=u_0(x),\notag\\
&&\lim_{|x|\to\infty}u(t,x)=0,\notag
\end{eqnarray}
can be given in terms of the vorticity field $\xi=\curl u$ as
\begin{eqnarray*}
&&\partial_t\xi-\nu\Delta\xi+(u\cdot\nabla)\xi-(\xi\cdot\nabla)u=g,\\
&&\xi(0,x)=\xi_0(x),\\
&&\xi=\curl u,\\
&&\Div u=0,\\
&&\lim_{|x|\to\infty}u(t,x)=0,
\end{eqnarray*}
where $g=\curl f$. We shall write the term $(\xi\cdot\nabla)u$ as
$(\nabla u)\xi$. Moreover, the same term can be written as $\D_u\xi$,
where $\D_u$ is the \textsl{deformation tensor}, the symmetric part of
$\nabla u$,
$$
\D_u=\frac12(\nabla u+\nabla u^T),
$$
since
$$
(\nabla u)\xi-\D_u\xi=\frac12(\nabla u-\nabla u^T)\xi=\xi\times\xi=0.
$$
As we explained intuitively in the introduction (see Section
\ref{interpreto}) and we shall describe rigorously in the sequel,
using the representation formula of Theorem \ref{biot} and the
generalised Feynman-Kac formula of Theorem \ref{fmain}, the formulation
of Navier-Stokes equations can be given in the following way:
\begin{eqnarray}\label{probform}
&&\xi(t,x)=\E[U^{x,t}_t\xi_0(X^{x,t}_t)]
          +\int_0^t\E[U^{x,t}_sg(t-s,X^{x,t}_s)]\,ds,\\
&&u(t,x)=\frac12\int_0^\infty\frac1s\E[\xi(t,x+W_s)\times W_s]\,ds,\notag
\end{eqnarray}
where the \textsl{Lagrangian paths} $(X^{x,t}_s)_{0\le s\le t}$ are
processes solutions of the following stochastic differential equations
$$
\left\{\begin{array}{ll}
dX^{x,t}_s=-u(t-s,X^{x,t}_s)\,ds+\sqn\,dW_s,&\qquad s\le t,\\
X^{x,t}_0=x,
\end{array}\right.
$$
and the \textsl{deformation matrices} $(U^{x,t}_s)_{0\le s\le t}$ are
the solutions to the following differential equations with random
coefficients
$$
\left\{\begin{array}{ll}
dU^{x,t}_s=U^{x,t}_s\D_u(t-s,X^{x,t}_s)\,ds,&\qquad s\le t,\\
U^{x,t}_0=I,
\end{array}\right.
$$
Here $\D_u$ is either $\nabla u$ or the deformation tensor (the name
deformation matrices of $U^{x,t}_s$ refers to the latter case). Notice
that, with respect to the introduction, we have made a time-reversion
which simplifies the mathematical analysis.

A sufficiently regular solution of the classical formulation
\eqref{nsformulation} is a solution of \eqref{probform} and
vice-versa. The main aim of this section is to show that, under suitable
conditions, problem \eqref{probform} has a unique local in time
solution. The claim is proved in the following theorem.

\begin{theorem}
Given $p\in[1,\frac32)$, $\alpha\in(0,1)$ and $T>0$, let $\xi_0\in
C^\alpha_b(\erre^3,\erre^3)\cap L^p(\erre^3,\erre^3)$ and
$g\in\valfapi$, and set
$$
\ep_0=\|\xi_0\|_{C^\alpha_b\cap L^p}+\int_0^T\|g(s)\|_{C^\alpha_b\cap L^p}\,ds.
$$
Then there exists $\tau\in(0,T]$, depending only on $\ep_0$, such that
there is a unique solution $u\in\tualfa$, with $\xi\in\tvalfapi$, of
problem \eqref{probform}.
\end{theorem}
\begin{proof}
We will show that there are suitable $L$, $M$ and $\tau$ such that the
map $\bs\circ\ns$, where $\bs:\tvalfapil\to\tualfam$ is defined as
$$
\bs(\xi)(t,x)=\frac12\int_0^\infty\frac1s\E[\xi(t,x+W_s)\times W_s]\,ds,
$$
and $\ns:\tualfam\to\tvalfapil$ is defined as
$$
\ns(u)(t,x)=\E[U^{x,t}_t\xi_0(X^{x,t}_t)]+\int_0^t\E[U^{x,t}_tg(t-s,X^{x,t}_s)]\,ds,
$$
is contractive.

First, in view of Corollary \ref{biotsavart}, $M\ge\widetilde CL$. Using
Proposition \ref{nsmap}, we see that $\ns$ maps $\tualfam$ to
$\tvalfapil$ if
\begin{equation}\label{conduno}
\e^{3\tau M}(1+\tau M)\ep_0\le L.
\end{equation}
By means of Corollary \ref{biotsavart} and Proposition \ref{nslip},
$\bs\circ\ns$ is contractive if
\begin{equation}\label{conddue}
\widetilde C C(\nu,p)C_M(\tau)\ep_0<1,
\end{equation}
where $C(\nu,p)$ is a constant depending only on $p$ and $\nu$, and
$\lim_{\tau\to 0}C_M(\tau)=0$. Hence, it is sufficient to choose
$\tau$ small enough in order to have both conditions \eqref{conduno} and
\eqref{conddue} verified.
\end{proof}

\begin{remark}
As usual, the statement of the above theorem can be read in terms of
small initial data. More precisely, for each fixed time $T$, there is a
constant $\epsilon$ such that if $\ep_0\le\epsilon$, there exists a
unique solution $u\in\ualfa$, with $\xi\in\valfapi$, of problem
\eqref{probform}
\end{remark}


\section{The Feynman-Kac formula for a deterministic system of parabolic
equations}\label{genpara}

This section is devoted to the development of a probabilistic
representation formula for the following system of parabolic equations
with final condition:
\begin{equation}\label{fgensys}
\left\{\begin{array}{l}
\partial_t v_k 
 +\frac12\sum_{i,j}a_{ij}\partial^2_{x_ix_j} v_k
 +\sum_{i=1}^d b_i\partial_{x_i} v_k
 +(\D v)_k
 +f_k=0,\\
v_k(T,x)=\varphi_k(x),\qquad x\in\erre^d,
\end{array}\right.
\quad k=1,\ldots,l,
\end{equation}
for $(t,x)\in[0,T]\times\erre^d$, or the following system of parabolic
equations with initial condition
\begin{equation}\label{igensys}
\left\{\begin{array}{l}
\partial_t v_k
= \frac12\sum_{i,j}a_{ij}\partial^2_{x_ix_j} v_k
 +\sum_{i=1}^d b_i\partial_{x_i} v_k
 +(\D v)_k
 +f_k,\\
v_k(0,x)=\varphi_k(x),\qquad x\in\erre^d,
\end{array}\right.
\quad k=1,\ldots,l,
\end{equation}
where $a=\sigma\sigma^*$ and
\begin{eqnarray}\label{data}
        \sigma:[0,T]\times\erre^d\longrightarrow\erre^{d\times d},
&\quad& b:[0,T]\times\erre^d\longrightarrow\erre^d,\notag\\
        \D:[0,T]\times\erre^d\longrightarrow\erre^{l \times l},
&     & \varphi:\erre^d\longrightarrow\erre^l\\
        f:[0,T]\times\erre^d\to\erre^l\notag
\end{eqnarray}
are Borel measurable functions. Additional assumptions will be stated in
the sequel.

At first, for simplicity, assume that $f\equiv0$ and all the data are
regular. If $l=1$, the equation \eqref{fgensys}, with final condition,
has a unique solution given by the Feynman-Kac formula
$$
v(t,x)=\E[\varphi(X_T^{t,x})\e^{\int_t^T\D(r,X_r^{t,x})\,dr}]
$$
where $X_s^{t,x}$ is the solution of the SDE
\begin{equation}\label{chara}
\begin{cases}
dX_s^{t,x}=b(s,X_s^{t,x})\,ds+\sigma(s,X_s^{t,x})\,dW_s, 
 &\quad s\in[t,T], \\
X_t^{t,x}=x,
\end{cases}
\end{equation}
where $(W_t)_{t\ge0}$ is a $d$-dimensional Brownian motion on some
filtered probability space. Our aim is to extend such formula to the
case $l>1$.

Notice that in the case $l=1$, for each $\omega$, the function
$$
u_r^{t,x}=\e^{\int_t^r\D(s,X_s^{t,x})\,ds},
$$
is the solution of the following equation (now $\D$ is a scalar)
\begin{equation}\label{1deformeq}
\begin{cases}
du_r^{t,x}=u_r^{t,x}\D(r,X_r^{t,x})\,dr,
 &\qquad r\in[t,T],\\
u_t^{t,x}=1.
\end{cases}
\end{equation}
So, in the same way, in the case $l>1$, we will consider the process
$U^{t,(x,Y)}$, solution of the equation
\begin{equation} \label{deformeq}
\begin{cases}
dU_r^{t,(x,Y)}=U_r^{t,(x,Y)}\D(r,X_r^{t,x})dr
 &\qquad r\in[t,T],\\
U_t^{t,(x,Y)}=Y,
\end{cases}
\end{equation}
where now both $\D$ and $U^{t,(x,Y)}$ are $l\times l$ matrices. If
$Y\equiv I$ we will write $U^{t,x}$ in place of $U^{t,(x,I)}$. Now, the
natural conjecture is that, under suitable regularity conditions, the
solution of \eqref{fgensys} is
\begin{equation} \label{lfk}
v(t,x)=\E[U^{t,x}_T \varphi (X_T^{t,x})].
\end{equation}

In Section \ref{shomo} we shall prove \eqref{lfk}, under suitable regularity
conditions on the coefficients. Such formula needs to be modified in
order to handle the case $f\not\equiv 0$, as we show in Section
\ref{sinhomo}. In Section \ref{suniq} we shall provide sufficient
conditions for the uniqueness of strong solutions to system
\eqref{fgensys}. Finally in Section \ref{sini} we shall give a
Feynman-Kac representation for the solutions of the system, with
initial condition, \eqref{igensys}.

\begin{remark}
When $l=1$, we can write without distinction in formula
\eqref{1deformeq} both $u_r^{t,x}\D$ and $\D u_r^{t,x}$, since they are
both scalars. If $l>1$, the lack of commutativity for the matrix
products gives that $U_r^{t,x}\D$ and $\D U_r^{t,x}$ are different. The
choice in the order of the matrix product in equation \eqref{deformeq},
and in formula \eqref{lfk} as well, derives from the form of the term
$\D\cdot v$ in system \eqref{fgensys}. To have an intuitive idea of
this fact, the reader can see the computations in the proof of the
uniqueness in Proposition \ref{stronguniq} (it is convenient to take
$f\equiv0$ for simplicity). However, when one uses backward stochastic
equations to represent solutions, the order of $U_r$ and $\D$ in
equation \eqref{deformeq} changes, see Section \ref{sini}.
\end{remark}


\subsection{The homogeneous case}\label{shomo}

Throughout this section, we will assume
$$
f\equiv 0$$
and that the functions $b$, $\sigma$ and $\D$, given in \eqref{data},
are Borel measurable functions such that
\begin{enumerate}
\item[\auno]{$b$, $\sigma$ are sub-linear with respect to $x$, uniformly in
$t$,}
\item[\adue]{$b$, $\sigma$ are locally Lipschitz-continuous in $x$,
uniformly in $t$,}
\item[\atre]{$a$ is differentiable in $x$ and $\partial_{x_i}a$ are locally
Lipschitz-continuous in $x$, uniformly in $t$,}
\item[\aqua]{$\D$ is bounded and locally Lipschitz-continuous in $x$,
uniformly in $t$,}
\item[\acin]{$\varphi$ is bounded and continuous.}
\end{enumerate}
In particular, assumptions \auno, \adue\ and \aqua\ ensure the existence
of strong solutions, unique in law, for the equations \eqref{chara} and
\eqref{deformeq}. Moreover, from assumption \aqua, it easily follows that
\begin{equation}\label{Ubound}
\|U_T^{t,x}\|_{\erre^{l\times l}}\le \e^{T\|\D\|_\infty},
\end{equation}
where $\|\D\|_\infty$ is the sup-norm. Finally, the previous formula and
assumption \acin, imply that the function $v$ given by formula
\eqref{lfk} is well defined and bounded.

We can now state the main result of this section:

\begin{theorem} \label{genmain}
Assume \auno-\acin\  and $\varphi\in C_b(\erre^d,\erre^l)$.
Then the function 
$$
v(t,x)=\E[U_T^{t,x}\varphi(X_T^{t,x})]
$$ 
is continuous and bounded and solves the Kolmogorov equation \eqref{fgensys}
in the sense of distributions, that is
\begin{equation} \label{distFK}
\int_0^T \int_{\erre^d}vM^*\eta\,dx\,dt=0,\qquad\textrm{for all }\eta\in
C_c^{\infty}((0,T)\times\erre^d,\erre^l),
\end{equation}
where
\begin{equation}\label{mop}
M^*\eta=- \partial_t\eta
          + \frac12\sum_{ij}\partial^2_{x_ix_j}(a_{ij}\eta)
          - \sum_i \partial_{x_i}(b_i \eta)
          + \D^*\eta.
\end{equation}
\end{theorem}

\begin{remark}\label{opsense}
The operator $M^*$ makes sense since, by assumptions \adue,
\atre\ and Rademacher theorem, the functions $\partial_{ij}a_{ij}$ and
$\partial_i b_i$ are well defined a.e.\ and essentially bounded in
compact sets. Moreover, $M^*\eta$ is bounded in compact sets.
\end{remark}

To prove Theorem \ref{genmain}, we shall use the  \textsl{method of new
variables} given by Krylov in \cite{Kry}.  Krylov used such method in
order to transform a parabolic equation on $\erre^d \times [0,T]$ with
potential term, into a parabolic equation on $\erre^{d+2} \times [0,T]$
without potential term. As observed in the introduction, we extend this
method to systems of parabolic equations. In our case, the elimination
of the potential term has the additional advantage that the coupling
between the equations in \eqref{fgensys} disappears. In other words,
we turn system \eqref{fgensys} into a system of $l$
\textsl{independent} parabolic equations on $\erre^{d+l\times l}\times
[0,T]$ without the potential term.

We define the new variables $\ox=(x,Y)\in\erre^{d+l\times l}$, and,
for each function $\psi:\erre^d\to\erre^l$, we define the function
$\overline{\psi}:\erre^{d+l\times l}\to\erre^l$ as
$\overline{\psi}(\ox)=Y \psi(x)$. Finally, if
$u(t,x):[0,T]\times\erre^d\to\erre^l$, we set $\ou(t,\ox)=Yu(t,x)$.

Prior to the computation of the derivatives of $\ov$, we give
some notations. We denote by $\zero{m}{n}$ the $m\times
n$ matrix with all entries equal to zero. Given a column vector
$\alpha\in\erre^d$ and a $l\times l$ matrix $A$, we define
the $(d+l)$ (\textsl{exotic}) column vector
$\left[\begin{smallmatrix}\alpha\\A\end{smallmatrix}\right]$,
where the first $d$ rows are given by the components of $\alpha$ and
the other $l$ rows are the rows of $A$ (the apparent inconsistency is
inessential, since we shall only use the scalar product defined below).
The scalar product between two such vectors is defined as
$$
\langle\left[\begin{smallmatrix}\alpha\\A\end{smallmatrix}\right],
\left[\begin{smallmatrix}\beta\\B\end{smallmatrix}\right]\rangle=
\alpha\cdot\beta+\langle A:B\rangle,
$$
where, as usual, $\langle A:B\rangle=\Tr(A\cdot B)=\sum_{i,j=1}^lA_{ij}B{ij}$.

Given $u\in C^1([0,T]\times\erre^d;\erre^l)$, since
\begin{equation}\label{Ygrad}
\frac{\partial \ou_h}{\partial Y_{ij}}
=\frac{\partial(Yu)_h}{\partial Y_{ij}}
=\frac{\partial}{\partial Y_{ij}}\sum_k Y_{hk}u_k=\delta_{ih}u_j,
\qquad h=1,\ldots,l
\end{equation}
it follows that, for each $h=1,\ldots,l$, the gradient $\nabla_\ox\ou_h$
of $\ou_h$ with respect to all its variables is given by the following
(exotic) column vector
$$
\nabla_\ox\ou_h
=\begin{bmatrix}
   \nabla_x (Yu)_h\\
   \nabla_{Y}\ou_h
 \end{bmatrix}
=\begin{bmatrix}
   \nabla_x (Yu)_h\\
   \zero1l\\\ldots\\
   u\\
   \ldots\\\zero1l
 \end{bmatrix},
$$
where the $d$-column vector is the gradient with respect to $x$ and the
$l\times l$ matrix has its rows all equal to the $l$-dimensional vector
$\zero1l=(0,\ldots,0)^T$ except for the $h^{\textrm{th}}$, which is the vector
$u$.

We want to evaluate next the scalar product
$\langle\left[\begin{smallmatrix}b\\Y\D\end{smallmatrix}\right],\nabla_\ox\ou_h\rangle$.
Since 
$$
(Y\D)_{ij}\partial_{Y_{ij}}(Yu)_h
=(Y\D)_{ij}\delta_{ih}u_j
=(Y\D)_{hj}u_j\delta_{ih},
$$
it follows that
$$
\langle\left[\begin{smallmatrix}
         b\\
         Y\D
       \end{smallmatrix}\right]
,\nabla_\ox\ou_h\rangle
=b\cdot\nabla_x (Yu)_h +(Y\D u)_h.
$$
In particular, if $Y=I$, the above quantity is equal to $b\cdot\nabla_x
u_h+(\D u)_h$.
 
Let
\begin{equation}\label{newdd}
\alpha(t,\ox)=\begin{pmatrix}
                a(t,x)&\zero{d}{l^2}\\
                \zero{l^2}{d}&\zero{l^2}{l^2}
              \end{pmatrix},
\qquad
\beta(t,\ox)=\begin{bmatrix}
                b(t,x)\\
                Y\D(t,x)
             \end{bmatrix},
\end{equation}
where we understand that $\alpha$ is defined in blocks, where each entry
is a matrix itself (notice that also $D^2_\ox\ou_h$ is defined in
blocks, and the product $\langle\alpha:D^2_\ox\ou_h\rangle$ is defined
as the sum of the four $\langle\,:\,\rangle$-products of the
corresponding blocks). With these notations, if $u\in
C^{1,2}([0,T]\times\erre^d;\erre^l)$, we have for each $h=1,\ldots,l$,
\begin{eqnarray*}
\lefteqn{\partial_t\ou_h+\frac12\langle\alpha:D^2_\ox\ou_h\rangle+\langle\beta,\nabla_\ox\ou_h\rangle=}\\
&=& (Y\partial_t u)_h
    +\frac12\sum_{i,j}a_{ij}\partial^2_{x_ix_j}(Yu)_h
    +\sum_i b_i\partial_{x_i}(Yu)_h
    +(Y\D u)_h\\
&=& \Big[Y\big(\partial_t u
                 +\frac12\sum_{i,j}a_{ij}\partial^2_{x_ix_j}u
                 +\sum_i b_i\partial_{x_i}u
                 +\D u\big)\Big]_h.
\end{eqnarray*}
From this identity it is straightforward to prove that a field $u$ is a
strong solution of \eqref{fgensys} if and only if $\ou$ is a strong
solution of system \eqref{newgen}, where by strong solution we mean a
continuous function having continuous first derivatives in time and
second derivatives in space, and satisfying the corresponding equation
point-wise. In the same way, applying the same ideas used above on the
adjoint operator, we have the following equivalence.

\begin{proposition}\label{weakweak}
A function $u$ is a weak solution of system \eqref{fgensys}, with
final condition, if and only if $\ou$ is a weak solution of
\begin{equation} \label{newgen}
 \partial_t\ou_h
+\frac12\langle\alpha:D^2_\ox\ou_h\rangle
+\langle\beta,\nabla_\ox\ou_h\rangle
=0,\qquad h=1,\ldots,l,
\end{equation}
with final condition $\ou(T,\ox)=Y\varphi(x)$.
\end{proposition}

In the sequel we prove that, under suitable conditions, the vector field
$\ov(t,\ox)=Yv(t,x)$, where $v$ is given by \eqref{lfk}, is a weak
solution of \eqref{newgen}. In view of the above lemma, this implies
that the function given by \eqref{lfk} solves system \eqref{fgensys} in
the weak sense.

The main part is contained in the following proposition, where we relax
some regularity assumptions on the coefficients of a theorem of
Krylov \cite{Kry}. Indeed, the drift and the diffusion defined in formulae
\eqref{newdd} are neither bounded nor globally Lipschitz-continuous, in
contrast to the assumptions of \cite{Kry}. The same problem occurs for
the final condition. On the other hand, both the drift and the
diffusion are locally Lipschitz-continuous and with linear growth (in
all variables, including $Y$).

\begin{proposition}\label{krylovrep}
Let $m\in\N$ and consider the scalar parabolic equation
\begin{equation} \label{krygen}
 \partial_tu
+\frac12\langle\alpha:D^2u\rangle
+\langle\beta,\nabla u\rangle
=0\qquad (t,x)\in[0,T]\times\erre^m
\end{equation}
with final condition $u(T,x)=\psi(x)$, where $\alpha=\gamma\gamma^*$ and
$$
\beta:[0,T]\times\erre^m\to\erre^m,
\quad
\gamma:[0,T]\times\erre^m\to\erre^{m\times m},
\quad
\psi:[0,T]\times\erre^m\to\erre,
$$
and assume that
\begin{enumerate}
\item[(i)]{$\beta$, $\gamma$ are Borel measurable, sub-linear and
locally Lipschitz continuous in $x$, uniformly in $t$,}
\item[(ii)]{$\psi$ is continuous and with polynomial growth,}
\item[(iii)]{$\gamma(t,\cdot)$ is continuously differentiable for each
$t$ and $\partial_{x_i}\gamma$ are locally Lipschitz continuous in x, uniformly
in $t$.}
\end{enumerate}
Set $u(t,x)=\E[\psi(Z^{t,x}_T)]$, where $Z^{t,x}_r$ is the solution of
the SDE
$$
\begin{cases}
dZ_r^{t,x}=\beta(r,Z_r^{t,x})\,dr + \gamma(r,Z_r^{t,x})\,dW_r,&\quad r\in[t,T],\\
Z_t^{t,x}=x,
\end{cases}
$$
where $(W_t)_{t\ge0}$ is an $m$-dimensional standard Brownian motion.
Then $u$ is a weak solution of \eqref{krygen}: for each $\eta\in
C_c^\infty((0,T)\times\erre^m)$, we have
$$
\int_0^T\int_{\erre^m}uN^*\eta\,dx\,dt=0,
$$
where
$$
N^*\eta=- \partial_t\eta
        + \frac12\sum_{i,j} \partial^2_{x_ix_j}(\alpha_{ij}\eta)
        - \sum_i \partial_{x_i}(\beta_i\eta).
$$
\end{proposition}

\begin{proof}
If everywhere in the assumptions of the proposition we have global
Lipschitz-continuity (instead of local Lipschitz-continuity), the
proposition follows from Theorem 5.13 of Krylov \cite{Kry}. In
the general case, we proceed by truncation. Let $\Psi_n\in
C^{\infty}(\erre^m)$ be such that
$$
\Psi_n(x)=\begin{cases}1&\quad |x|\le n\\
                       0&\quad |x|\ge n+1
          \end{cases}
$$
and set $\beta^{(n)}=\Psi_n\beta$ and $\gamma^{(n)}=\Psi_n\gamma$.
Fix a Brownian motion $(\Omega,\F,\F_t,W_t,\Pb)$ and
denote by $Z_t^{s,x,n}$ the solutions to the corresponding SDEs.
The sequence $Z_t^{s,x,n}$ converges to $Z_t^{s,x}$ in probability
uniformly on compact subsets of $[0,T]\times\erre^m$.

Suppose first that $\psi$ is bounded. Then
$u_n(t,x)=\E[\psi(Z_t^{s,x,n})]$ converges to $u(t,x)$ and $\beta_{x_i}^{(n)}$
converges to $\beta_{x_i}$, $\partial_{x_i}\alpha^{(n)}$ to
$\partial_{x_i}\alpha$ and $\partial_{x_i,x_j}\alpha^{(n)}$ to
$\partial_{x_i,x_j}\alpha$ uniformly on compact subsets of
$[0,T]\times\erre^m$. Let $\eta\in C_c^\infty$, since $N_n^*\eta$ is a
bounded sequence (see Remark \ref{opsense}), by the dominated
convergence theorem, $\int u_nN^*_n\eta$ converges to $\int uN^*\eta$,
where $N_n^*$ is the operator corresponding to the approximate
coefficients. Since $u_n$ are weak solutions, it follows that $u$ is
also a weak solution.

If $\psi$ is not bounded, we take a sequence of bounded continuous
functions $\psi_n\to\psi$ such that $|\psi_n(x)| \le |\psi(x)|$. From
Theorem 4.6 of Krylov \cite{Kry}, we have $\E[|Z_T^{t,x}|^k]\le
c(1+|x|^k)$, so that $u_n(t,x)\le c(1 +|x|^k)$ by assumption $(ii)$, and
again we conclude by the dominated convergence theorem.
\end{proof}

We are now ready to prove the main theorem.

\begin{proof}[Proof of Theorem \ref{genmain}]
First we show that $v$ is bounded and continuous. The boundedness comes
from \eqref{Ubound} and the assumptions on $\varphi$. In order to show
the continuity, we take a sequence $(x_n,t_n)$ converging to
$(x,t)$. From Lemma 2.9 of Krylov \cite{Kry}, the function
$(t,x)\to(X_\cdot^{t,x},U_\cdot^{t,(x,I)})\in C([0,T],\erre^{d+l\times
l})$ (where by convention $(X_s^{t,x}, U_s^{t,(x,I)})=(x,I)$ if $s<t$)
is continuous in probability. Hence, there is a subsequence such that
convergence is almost sure. Finally, the conclusion follows from the
bound \eqref{Ubound}, the assumptions on $\varphi$ and the dominated
convergence theorem.

We show then that $v$ is a weak solution. We have the following two
ingredients:
\begin{enumerate}
\item[{\it (i)}] the two systems of SDEs \eqref{chara} and
\eqref{deformeq} can be thought as a unique system where the solution
$(X_r^{t,x},U^{t,(x,Y)}_r)$ takes values in $\erre^{d+l\times l}$ and
drift and diffusion are given by formulae \eqref{newdd}.
\item[{\it (ii)}] Since by uniqueness for equation \eqref{deformeq} it
follows that $U_T^{t,(x,Y)}=YU_T^{t,x}$, for the function $v$ defined in
\eqref{lfk}, we have
$$
\ov(t,\ox)
= Yv(t,x)=\E[YU^{t,x}_T\varphi(X^{t,x}_T)]
= \E[U_T^{t,(x,Y)}\varphi(X_T^{t,x})]
=\E[\overline{\varphi}(X_T^{t,x},U_T^{t,(x,Y)})].
$$
\end{enumerate}
From these two facts, by Proposition \ref{krylovrep}, $\ov$ is a weak
solution to system \eqref{newgen}. By Proposition \ref{weakweak}, $v$ is
a weak solution to system \eqref{fgensys}.
\end{proof}

The regularity assumption \aqua\ on the term $\D$ can be relaxed with
the following condition
\begin{enumerate}
\item[\aquap]{$\D$ is bounded and uniformly continuous.}
\end{enumerate}
In fact we can deduce the following corollary.

\begin{corollary} \label{gencor}
Assume \auno-\atre, \aquap\ and \acin. Then the function
$$
v(t,x)=\E[U_T^{t,x}\varphi(X_T^{t,x})]
$$
is continuous and bounded and solves the Kolmogorov equation
\eqref{fgensys} in the sense of distributions:
$$
\int_0^T\int_{\erre^d}vM^*\eta\,dx\,dt=0
\qquad\text{for all }\eta\in C_c^\infty((0,T)\times\erre^d,\erre^d).
$$
\end{corollary}

\begin{proof}
Let $\rho_n$ be a sequence of mollifiers and set $\D_n=\D*\rho_n$ and
$v_n(t,x)=\E[U_{T,n}^{t,x}\varphi(X_T^{t,x})]$, where $U_{r,n}^{t,x}$ is
the solution of \eqref{deformeq} corresponding to $\D_n$.

Since $\D_n\to\D$ uniformly in $[0,T]\times\erre^d$,
we have $U_{t,n}^{t,x}\to U_t^{t,x}$ in $L^1(\Omega)$, uniformly in
$[0,T]\times\erre^d$. Consequently, $v_n(t,x)\to v(t,x)$ and $\D_n v_n
\to \D v$ uniformly $[0,T]\times\erre^d$. Since $v_n$ are weak solutions
of the corresponding approximate problem, in the limit $v$ is a weak
solution of $Mv=0$.
\end{proof}


\subsection{The inhomogeneous case}\label{sinhomo}

In this section, Theorem \ref{genmain} will be extended to the
inhomogeneous case. We will show a Feynman-Kac representation formula
for the complete system \eqref{fgensys}, that is with $f\not\equiv 0$,
with final condition. Throughout this section we will assume
\auno-\atre, \aquap, \acin\ and the following
\begin{enumerate}
\item[\asei]{$f:[0,T]\times\erre^d\to\erre^l$ is bounded and uniformly
continuous.}
\end{enumerate}
\begin{theorem} \label{inhomo}
Assume \auno-\atre, \aquap, \acin-\asei. Then the function
\begin{equation}\label{gfk}
v(t,x)=  \E[U_T^{t,x}\varphi(X_T^{t,x})]
       + \int_t^T\E[U_r^{t,x}f(r,X_r^{t,x})]\,dr
\end{equation}
is a weak solution of \eqref{genmain}, that is,
$$
\int_0^T\!\int_{\erre^d}\bigl(u M^*\eta+f\eta\bigr)\,dt\,dx=0,
\qquad\eta\in C_c^\infty((0,T)\times\erre^d,\erre^l).
$$ 
\end{theorem}

The main idea to prove the theorem is to introduce a new component (we
apply again the \textsl{method of new variables} of Krylov \cite{Kry})
and prove that $v$ is a solution of system \eqref{fgensys} if and only
if $\tv=(v_1,\ldots,v_l,1)$ solves the system
\begin{equation}\label{newinho}
 \partial_t\tv
+\frac12\sum_{i,j}a_{ij}\partial_{x_ix_j}\tv
+\sum_i b_i\partial_{x_i}\tv
+(\tD\tv)_k
=0,
\end{equation}
with final condition $\tv(T,\cdot)=(\varphi_1,\ldots,\varphi_l,1)$, where
$\tD=\big(\begin{smallmatrix}\D&f\\0&0\end{smallmatrix}\big)$. Notice that
$\tD\tv=\big(\begin{smallmatrix}\D v+f\\0\end{smallmatrix}\big)$, so that
the component $\tv_{l+1}$ is obviously a solution. 

The key lemma is the following.

\begin{lemma}
The function $\tilde v=(v_1,v_2,...,v_l,1)$ is a weak solution of
 \eqref{newinho} if and only if $v=(v_1,v_2,...,v_l)$
is a weak solution of \eqref{fgensys}.
\end{lemma}

\begin{proof}
A weak solution of \eqref{fgensys} is a function $v$ such that
$\iint (vM^*\eta+f\eta)=0$ for each test function $\eta$, or equivalently
$\iint (vL^*\eta+v\D^*\eta+f\eta)=0$, where the operator $M^*$ has been
defined in \eqref{mop} and $L^*$ is defined as
$$
L^*\eta=-\partial_t\eta
        +\frac12\sum_{i,j}\partial^2_{x_ix_j}(a_{ij}\eta)
        -\sum_i\partial_{x_i}(b_i\eta).
$$
Let $\tilde\eta=(\eta,\eta_{l+1})$ be a $\erre^{l+1}$-valued test function.
Since $\D^*=\big(\begin{smallmatrix}\D^*&0\\f^*&0\end{smallmatrix}\big)$,
we have
$$
\tv\tD^*\tilde\eta=\begin{pmatrix}
                  v\D^*\eta+f\eta\\
                  0
                 \end{pmatrix}.
$$
It comes out that $v$ is a solution of the inhomogeneous equation if and
only if $\tv$ solves $\iint(\tv L^*\tilde\eta+f\tilde\eta)=0$, that is, if
and only if $\tv$ is a weak solution of system \eqref{newinho}.
\end{proof}

We can now prove the main theorem of this section.

\begin{proof}[Proof of Theorem \ref{inhomo}]
Let $\tilde\varphi$ be the function $(\varphi_1,...,\varphi_l,1)$ and
$\tilde U_s^{t,x}$ be the solution of 
\begin{equation} \label{indeformeq}
\begin{cases}
d\tU_s^{t,x}=\tU_s^{t,x}\tD(s,X_s^{t,x})\,ds, &\qquad s\in[t,T], \\
\tU_t^{t,x}=I_{l+1}.
\end{cases}
\end{equation}
Since $\varphi$, $\D$ and $f$ satisfy assumptions \aquap, \acin\ and
\asei, the functions $\tilde\varphi$ and $\tD$ satisfy assumptions
\aquap\ and \acin. Hence, by Corollary \ref{gencor}, the function 
$$
(x,t)\to\E[\tU_T^{t,x}\tilde\varphi(X_T^{t,x})]
$$ 
is a weak solution of system \eqref{newinho}.

We write $\tU^{t,x}_s$ in blocks:
$$
\tU_s^{t,x}=\begin{pmatrix}
              A_s^{t,x}&b_s^{t,x}\\
              c_s^{t,x}&d_s^{t,x}
            \end{pmatrix},
$$
where $A_s$ is a $l\times l$ matrix, $b_s\in\erre^d$ is a column vector,
$c_s\in\erre^d$ is a row vector and $d_s$ is a scalar. With this
position, the Cauchy problem (\ref{indeformeq}) is equivalent to
$$
\begin{cases}
dA_s^{t,x}=A_s^{t,x}\D(s,X_s^{t,x})\,ds, &\qquad A_t^{t,x}=I_{l},\\
db_s^{t,x}=A_s^{t,x}f(s,X_s^{t,x})\,ds,  &\qquad b_t^{t,x}=0,\\
dc_s^{t,x}=c_s^{t,x}\D(s,X_s^{t,x})\,ds, &\qquad c_t^{t,x}=0,\\
dd_s^{t,x}=c_s^{t,x}f(s,X_s^{t,x})\,ds,  &\qquad d_t^{t,x}=1,
\end{cases}
$$
and it is easy to see that
\begin{eqnarray*}
&A_s^{t,x}=U_s^{t,x} &\qquad b_s^{t,x}=\int_t^sU_r^{t,x}f(r,X_r^{t,x})\,dr\\
&c_s^{t,x}=0&\qquad d_s^{t,x}=1.
\end{eqnarray*}
Consequently,
\begin{eqnarray*}
\E[\tU_T^{t,x}\tilde\varphi(X_T^{t,x})]
&=&\E\big[\big(\begin{smallmatrix}U_T^{t,x}&b_T^{t,x}\\0&1\end{smallmatrix}\big)
      \big(\begin{smallmatrix}\varphi(X_T^{t,x})\\1\end{smallmatrix}\big)\big]
=\E\big[\big(\begin{smallmatrix}U^{t,x}_T\varphi(X_T^{t,x}+b^{t,x}_T)\\1\end{smallmatrix}\big)\big]\\
&=&\E\big[\big(\begin{smallmatrix}U^{t,x}_T\varphi(X_T^{t,x}+b^{t,x}_T)+\int_t^sU_r^{t,x}f(r,X_r^{t,x})\,dr\\1\end{smallmatrix}\big)\big].
\end{eqnarray*}
\end{proof}


\subsection{A uniqueness result}\label{suniq}

In the preceding sections, we were concerned with the existence of a
weak solution of the parabolic system \eqref{fgensys} having a nice
probabilistic representation. The aim of the present section is to
provide sufficient conditions for the uniqueness of solutions. In
Proposition \ref{stronguniq} we shall see that the strong solution, if
exists, is given by our probabilistic representation, hence is unique.
In Theorem \ref{weakuniq} we will show, under some special conditions on
the coefficients, that weak solutions are also unique and are given by
the probabilistic representation. Such special conditions on the
coefficients are satisfied in the application of the probabilistic
representation to the Navier-Stokes system: if the velocity field is
regular enough, the coefficients in the equations for the vorticity
satisfy the special conditions. Hence, for each fixed regular velocity,
there exists a unique weak solution of the vorticity equation given
by the Feynman-Kac formula.
\medskip

Let $C^{1,2}_b([0,T]\times\erre^d,\erre^l)$ be the space of continuous
functions having first and second derivatives in $x$ and first
derivative in $t$ continuous and bounded. We start by showing that, if
the solution of the parabolic system is regular, then it is given by
formula \eqref{gfk}.

\begin{proposition}\label{stronguniq}
Let $v\in C^{1,2}_b([0,T]\times\erre^d,\erre^l)$ be a strong solution of
system \eqref{fgensys}, with final condition. Then $v$ is given by
formula \eqref{gfk}.
\end{proposition}
\begin{proof}
It is sufficient to show that the process
$$
U_r^{t,x}v(r,X_r^{t,x})+\int_t^rU^{t,x}_sf(s,X^{t,x}_s)\,ds,
\qquad r\in[t,T],
$$
is a martingale. Indeed, if  $h\in\{1,\ldots,l\}$, by
It\^o formula (we omit for simplicity $(r,X_r^{t,x})$ from the term
$v(r,X_r^{t,x})$ and from the coefficients, and the subscript $r$ from
the term $U_r^{t,x}$)
\begin{eqnarray*}
d_r(U^{t,x}v)_h
&=& \sum_k d(U_{hk}^{t,x}v_k)
 =  \sum_k (U_{hk}^{t,x}dv_k
    + v_kdU_{hk}^{t,x})\\
&=& \sum_k U_{hk}^{t,x}\big[(\partial_r v_k
    + \sum_i b_i\partial_{x_i} v_k
    + \frac12\sum_{i,j}a_{ij}\partial^2_{x_ix_j}v_k)\,dr\\
& & + \sum_{i,j} \partial_{x_i} v_k\, \sigma_{ij}\,dW^j_r\big]
    +\sum_{k,i} v_kU^{t,x}_{hi}\D_{ik}\,dr\\
&=& -\sum_k U_{hk}^{t,x}\big(f_k+\sum_i\D_{ki}v_i\big)\,dr
    +(dM_r)_h
    +\sum_{k,i} v_kU^{t,x}_{hi}\D_{ik}\,dr\\
&=& -d_r\big(\int_t^r(U^{t,x}_sf)_h\,ds\big)
    +(dM_r)_h
\end{eqnarray*}
since $v$ is a solution of system \eqref{fgensys}; $(M_r)_{r\in[t,T]}$
is the $d$-dimensional martingale, vanishing at $r=t$, given by
$$
(dM_r)_h=\sum_k U^{t,x}_{hk}\sum_{i,j} \partial_{x_i}v_k\, \sigma_{ij}\,dW^j_r.
$$
Moreover, $M_r$ is square-integrable, since $v\in C^{1,2}_b$,
 $U_T^{t,x}$ is bounded by \eqref{Ubound} and
$$
\sup_{t\le r\le T}\E[|X_r^{t,x}|^2]
$$
is bounded.
\end{proof}
\begin{theorem}\label{weakuniq}
Let $\varphi$ be bounded and continuous, $f$ and $\D$ bounded and
uniformly continuous. Suppose that $\sigma$ is constant and $b$ is a Borel
measurable and Lipschitz-continuous in $x$ function such that $\Div b=0$.
Then the function
$$
v(t,x)
= \E[U_T^{t,x}\varphi(X_T^{t,x})]
 +\int_t^T\E[U_r^{t,x} f(r,X_r^{t,x})]\,dr
$$
is the unique weak solution of the parabolic system \eqref{fgensys}.
\end{theorem}
The proof of the theorem is based on a regularisation by convolution, in
order to apply the uniqueness result of the previous proposition.

Let $\rho\in C^\infty(\erre^d,\erre)$, $0\le\rho\le1$, with support in
the ball of radius one, and set $\rho_n(x)=n^d\rho(nx)$. Let $J_n$ be
the convolution operator: $J_n(u)=\rho_n*u$.
 
\begin{lemma}
Let  $b:\erre^d\to\erre^d$ be a Lipschitz-continuous function, such that
$\Div b=0$ (in the sense of distributions). Then there is a constant $C$
such that for each $u\in C_b(\erre^d,\erre^l)$,
\begin{equation}\label{urel}
\big|\big([J_n,b\cdot\nabla]u\big)(x)\big|\le C\sup_{y\in B_{1/n}(x)}
|u(y)|\qquad\text{for all }n,
\end{equation}
where $[J_n,b\cdot\nabla]u=J_n((b\cdot\nabla)u)-(b\cdot\nabla)J_nu$ is the commutator. Moreover
\begin{equation}\label{drel}
[J_n,b\cdot \nabla]u\stackrel{n\to\infty}{\longrightarrow}0
\qquad\text{uniformly on compact sets}
\end{equation}
\end{lemma}
\begin{proof}
Fix $u\in C_b(\erre^d,\erre^l)$. Since $\Div b=0$, by integration by
parts we have
\begin{eqnarray*}
\lefteqn{\big([J_n,b\cdot\nabla]u\big))(x)=}\\
&=& \big(\rho_n*(b\cdot\nabla)u-(b\cdot\nabla)(\rho_n*u)\big))(x)\\
&=& \int_{\erre^d}\rho_n(x-y)(b(y)\cdot\nabla_y)u(y)
     -(b(x)\cdot\nabla_x)(\rho_n(x-y))u(y)\,dy\\
&=& \int_{\erre^d}\nabla_y\rho_n(x-y)(b(x)-b(y))u(y)\,dy.
\end{eqnarray*}
Taking the norms in $\erre^l$ we get
\begin{eqnarray*}
\big|([J_n,b\cdot\nabla]u)(x)\big|
&\le& \int_{B_{1/n}(x)}|\nabla\rho_n(x-y)|\cdot|b(y)-b(x)|\cdot|u(y)|\,dy\\
&\le& cL\|\nabla\rho\|_\infty\sup_{y\in B_{1/n}(x)}|u(y)|
\end{eqnarray*}
where $L$ is the Lipschitz constant of $b$. So far, we have proved
\eqref{urel}. Concerning \eqref{drel}, it is easy to see that the
claim is true for $u\in C_b^\infty(\erre^d,\erre^l)$. If $u$ is only
$C_b$, the claim follows from approximation with $C^\infty_b$ functions
(in the sup-norm, on compact sets) and from the bound \eqref{urel}.
\end{proof}

We apply now the previous lemma to prove the main theorem.

\begin{proof}[Proof of Theorem \ref{weakuniq}]
Let $v$ be a bounded and continuous weak solution of system \eqref{fgensys}.
The sequence $v_n=\rho_n*v$ belongs to
$C([0,T],C_b^\infty(\erre^d,\erre^l))$ and $v_n\to v$ uniformly on
compact sets. We want to show that $v_n$ is a weak solution of
\begin{equation} \label{genapprox}
 \partial_t v_n
+\frac12\sum_{i,j}a_{ij}\partial^2_{x_ix_j} v_n
+\sum_ib_i\partial_{x_i} v_n
+\D v_n+\rho_n*f+R_n=0,
\end{equation}
with final condition $v_n(T)=\rho_n*\varphi$, where
$R_n=[J_n,b\cdot\nabla]v+[J_n,\D]v$. Indeed, $v$ is a weak solution of
\eqref{fgensys}, so that we can use $\zeta_n=\breve\rho_n*\eta$ as a
test function, where $\eta$ is again a test function and
$\breve\rho_n(x)=\rho_n(-x)$, to obtain with some easy computations
\begin{eqnarray*}
\lefteqn{\int_0^T\!\int_{\erre^d}(vM^*\zeta+f\zeta)=}\\
&=& \iint v\big(-\partial_t\zeta_n
                +\frac12\sum_{i,j} a_{x_ix_j}\partial^2_{ij}\zeta_n
                -\sum_i b_i\partial_{x_i}\zeta_n
                +\D^*\zeta_n\big)
                +f\zeta_n\\
&=& \iint \bigl[v_n\big(-\partial_t\eta
                  +\frac12\sum_{i,j} a_{ij}\partial^2_{x_ix_j}\eta
                  -\sum_i b_i\partial_{x_i}\eta
                  +\D^*\eta\big)\bigr]\\
& & +\iint\eta\big(J_nf
                   +[J_n,b\cdot\nabla]v
                   +[J_n,\D]v\big)
\end{eqnarray*}
(notice that for each $u$, $\int u(\breve\rho_n*\eta)=\int\eta(\rho_n*u)$).

Since $v_n$ belongs to $C([0,T],C_b^\infty(\erre^d,\erre^l))$
and $\rho_n*f+R_n$ is bounded and continuous, we argue that the 
distributional derivative $\partial_t v_n$ is bounded
and continuous and, therefore, a strong derivative.
Hence $v_n\in C^{1,2}_b$ and it is a strong solution of
(\ref{genapprox}). Proposition \ref{stronguniq} yields
$$
v_n(t,x)= \E[U_T^{t,x}\rho_n*\varphi(X_T^{t,x})]
         +\int_t^T\E[U_r^{t,x}(\rho_n*f+R_n)(X_r^{t,x})]\,dr.
$$
It is easy to check that $\|[J_n,D]v\|_\infty\le2\|D\|_\infty\|v\|_\infty$ and
$[J_n,D]v\to 0$, uniformly on compact sets. Hence, by the previous
lemma, $R_n$ is bounded, independently of $n$, and $R_n\to0$ uniformly
on compact sets. Using \eqref{Ubound} and the dominated convergence
theorem, we obtain
$$
v(t,x)=\lim_{n\to\infty}v_n(t,x)
      =\E[U_T^{t,x}\varphi(X_T^{t,x})]+
       \int_t^T\E[U_r^{t,x} f(r,X_r^{t,x})]\,dr.
$$
\end{proof}


\subsection{The formula for parabolic systems with initial
condition}\label{sini} 

In this section we describe the probabilistic representation of weak
solutions to the system \eqref{igensys}, with initial condition. Indeed,
in the sequel we will use the results of this sections to give a
probabilistic representation for the solutions to the Navier-Stokes
equations, which is a parabolic equation with initial condition.

We will obtain the representation formula for the forward parabolic system
using the representation for the backward parabolic system and a time
inversion of the coefficients. To this aim, we will consider the following
stochastic differential equations
\begin{equation}\label{ichara}
\begin{cases}
dX_r^{s,x,t}=b(t-r,X_r^{s,x,t})\,dr+\sigma(t-r,X_r^{s,x,t})\,dW_r,
 &\quad r\in[s,t],\\
X_s^{s,x,t}=x,
\end{cases}
\end{equation}
and
\begin{equation}\label{ideformeq}
\begin{cases}
dU^{s,(x,Y),t}_r=U^{s,(x,Y),t}_r\D(t-r,X_r^{s,x,t})\,dr,
 &\quad r\in[s,t],\\
U^{s,(x,Y),t}_s=Y,
\end{cases}
\end{equation}
where, as usual, $U^{s,(x,Y),t}=U^{s,x,t}$ when $Y=I$.

\begin{theorem} \label{fmain}
Let the data $b$, $\sigma$, $\varphi$, $\D$ and $f$ satisfy assumptions
\auno-\atre, \aquap\ (in page \pageref{gencor}, \acin\ (in page
\pageref{shomo}) and \asei (in page \pageref{inhomo}). Then the
function
\begin{equation} \label{ffk}
v(t,x)=\E[U^{0,x,t}_t\varphi(X_t^{0,x,t})]
       +\int_0^t\E[  U^{0,x,t}_r f(t-r,X_r^{0,x,t})]\,dr
\end{equation}
is a weak solution of \eqref{igensys}, with initial condition.

Moreover, if $\sigma$ is constant and $b$ is globally 
Lipschitz-continuous in $x$, then $v$ is the unique weak solution.
\end{theorem}

\begin{proof}
Let $\tv(t,x)=v(T-t,x)$. If $v$ is a weak solution of \eqref{igensys},
by easy computations it follows that $\tv$ is a weak solution of
\begin{equation} \label{igensystem}
\begin{split}
 \partial_t\tv(t,x)
&+\frac12\sum_{i,j}a_{ij}(T-t,x)\partial^2_{x_ix_j}\tv(t,x)
+\sum_i b_i(T-t,x)\partial_{x_i}\tv\\
&+\D(T-t,x)\tv(t,x)
+f(T-t,x)=0,
\end{split}
\end{equation}
for $t\in[0,T]$, with final condition $\tv(T,x)=\varphi(x)$ (and vice-versa).

By Theorem \ref{inhomo}, a solution $\tv$ of \eqref{igensystem} is given
by the following formula
$$
\tv(t,x)= \E[U^{t,x,T}_T\varphi(X_T^{t,x,T})]
         +\int_t^T\E[U^{t,x,T}_r f(T-r,X_r^{t,x,T})]\,dr,
$$
where $U^{t,x,T}_r$ and $X_r^{t,x,T}$ are given respectively in
\eqref{ideformeq} and \eqref{ichara}. We can conclude that a solution
$v$ of the forward parabolic equation \eqref{igensys} is given by
$$
v(t,x)= \E[U^{T-t,x,T}_T\varphi(X_T^{T-t,x,T})]
       +\int_{T-t}^T\E[U^{T-t,x,T}_r f(T-r,X_r^{T-t,x,T})]\,dr
$$
Finally, one can easily check that, for each $r\in[T-t,T]$, the joint
law of the random variables $U^{T-t,x,T}_r$ and $X_r^{T-t,x,T}$ is
equal to the joint law of the random variables $U^{0,x,t}_{r+t-T}$ and
$X_{r+t-T}^{0,x,t}$. In conclusion, formula \eqref{ffk} holds.
\end{proof}

The representation formula above appears more complicated 
than the formula for parabolic systems with final condition \eqref{gfk}:
the stochastic processes $X_r$ in \eqref{gfk} are the solutions
of a fixed SDE corresponding to different initial conditions, while
the stochastic processes $X_r^{0,x,t}$ and $U^{0,x,t}_r$ in \eqref{ffk} solve
for each $t$ a different SDE. A different representation can be given,
which is more appealing at the heuristic level, even if less suitable
for stochastic calculus. 

Consider the following backward SDE
\begin{equation}\label{bchara}
Y^{t,x}_r=x+\int_r^tb(s,Y^{t,x}_s)\,ds + \int_r^t \sigma(s,Y^{t,x}_s)\,\hat d W_s,
\qquad r \in [0,t],
\end{equation}
where $\hat dW_s$ denotes the backward stochastic integral with respect to
the Brownian motion $W_s$ (see Kunita \cite{Kun} for the definition of the
backward integral). Notice that the final condition $Y^{t,x}_t=x$ has
been imposed here. Let $V^{s,t,x}_r$, $0\le s\le r\le t$ be the
solution of
\begin{equation}\label{bdeformeq}
\begin{cases}
dV^{s,t,x}_r=\D(r,Y^{t,x}_r)V^{s,t,x}_r\,dr,
 &\qquad r\in[s,t],\\
V^{s,t,x}_s=I
\end{cases}
\end{equation}

\begin{theorem}\label{backmain}
Under the same assumptions of the previous theorem, a weak solution of
the parabolic system \eqref{igensys}, with initial condition,
is given by the following formula
\begin{eqnarray}\label{bfk}
v(t,x)= \E[V^{0,t,x}_t\varphi(Y_0^{t,x})]
       +\int_0^t\E[V^{r,t,x}_t f(r,Y_r^{t,x})]\,dr,
\end{eqnarray}
where $Y^{t,x}$ and $V^{r,t,x}$ are given respectively by \eqref{bchara}
and \eqref{bdeformeq}.
\end{theorem}

\begin{remark}
We want to give an interpretation of the representation formula given
above. Suppose for clarity that $f\equiv 0$. Consider the trajectory
$Y^{t,x}_r(\omega)$ of a virtual particle which is in $x$ at time $t$,
transported by a velocity field and subject to a diffusion, and evaluate
$v(0,Y_0^{t,x}(\omega))=\varphi(Y_0^{t,x}(\omega))$. Then we take into
account, through the vector field $V^{0,t,x}_t$, the effects of the
tensor $\D$ along the given trajectory in the time interval
$[0,t]$. Finally, by taking the expectation, we consider the mean effect
of all virtual particles.
\end{remark}

Before giving the proof of the theorem, we need the following simple
lemma for the time inversion of a stochastic integral.

\begin{lemma}\label{newbm}
Let $(W_s)_{s\ge0}$ be a Brownian motion. Fix $t>0$ and set
$$
B_s=W_t -W_{t-s} \qquad s\in[0,t].
$$
Let $\F^W_s=\sigma\left(W_r\,|\,r\in[0,s]\right)$ and
$\F^B_{s,t}=\sigma\left(B_u-B_v\,|\,s\le v\le u\le t\right)$ and let
$g(s)$ be a continuous and bounded process adapted to the filtration
$F^W_s$. Then the process $f(s)=g(t-s)$, $s\in[0,t]$ is
$\F_{s,t}^B$-adapted and for all $a,b$ such that $0\le a\le b\le t$,
$$
\int_a^b g(s)\,d W_s=\int_{t-b}^{t-a}f(s)\,\hat dB_s.
$$
\end{lemma}

\begin{proof}
Since $B_u-B_v=W_{t-v}-W_{t-u}$, we have $\F_{t-s}^W=\F_{s,t}^B$ and
this gives the first statement.

Take now a sequence of partitions of the interval $[a,b]$: 
$$
\pi_n:\{a=s_0^n\le s_1^n\le\ldots\le s_{k_n}^n=b\}
$$ 
such that $|\pi_n|\to 0$. We have
\begin{eqnarray*}
\int_a^b g(s)\,dW_s
&=&\lim_{n\to\infty}\sum g(s_i^n)(W_{s_{i+1}^n}-W_{s_i^n})\\
&=&\lim_{n\to\infty}\sum g(t-r_i^n)(W_{t-r_{i+1}^n}-W_{t-r_i^n})\\
&=&\lim_{n\to\infty}\sum f(r_i^n)(B_{r_i^n}-B_{r_{i+1}^n})\\
&=&\int_{t-b}^{t-a}f(s)\,\hat dB_s,
\end{eqnarray*}
where $r_i^n=t-s_i^n$, $i=1\ldots k_n$.
\end{proof}

\begin{proof}[Proof of Theorem \ref{backmain}]
We need only to show that
$$
X^{0,x,t}_{t-r}=Y^{t,x}_r
\qquad\text{and}\qquad
U^{0,x,t}_{t-r}=V^{r,x,t}_t,\quad \Pb-\text{a.s.}
$$
since such formulas, formula \eqref{ffk} and a change of variables give
us \eqref{bfk}.

We prove the first equality. Fix a Brownian motion $(W_r)_{r\ge 0}$ and
consider the solution $X^{0,x,t}_r$ of equation \eqref{ichara}. By
Lemma \eqref{newbm} above, it follows that $X^{0,x,t}_{t-r}$ satisfies
the backward SDE \eqref{bchara} with respect to the Brownian motion
$B_s$ defined in Lemma \ref{newbm}. Since equation \eqref{bchara} has a
unique strong solution, we have the first equality.

We proceed to prove the second equality. Fix $\omega$ so that $r\to
 Y_r^{t,x}(\omega)$ is continuous. The key observation is that
$$
V_r^{s,t,x}(\omega)=V_r^{0,t,x}(\omega)(V_s^{0,t,x}(\omega))^{-1},
\qquad 0\le s\le r\le t,
$$
and it is true since
$$
d(V^{0,t,x}_r(\omega))^{-1}=-(V^{0,t,x}_r(\omega))^{-1}\D(r,Y_r^{t,x}(\omega)),
$$
with initial condition $(V^{0,t,x}_0(\omega))^{-1}=I$, so that it easy to
check that $V_r^{0,t,x}(\omega)(V_s^{0,t,x}(\omega))^{-1}$
satisfy equation \eqref{bdeformeq}. Finally, by evaluating
$$
d_rV_t^{r,t,x}(\omega)
=d_r\big[V_t^{0,t,x}(\omega)(V_r^{0,t,x}(\omega))^{-1}\big],
$$
we see that both $V_t^{r,t,x}(\omega)$ and $r\to U^{0,x,t}_{t-r}(\omega)$ solves
the ODE: 
$$
dU_r=-U_r\D(r,Y_r^{t,x}(\omega))\,dr\qquad r\in[0,t],
$$
with final condition $U_t=I$.
\end{proof}


\section{A probabilistic representation for the Newtonian potential and
the Biot-Savart law}\label{sbiotsavart}

In the present section we aim to give a probabilistic representation for
the velocity field of an incompressible fluid in terms of the vorticity
field $\xi=\curl u$.

Under suitable assumptions on $\xi$, the Poisson equation
$-\Delta\psi=\xi$ has a solution, given by
$$
\psi(x)=\frac1{4\pi}\int_{\erre^3}\frac{\xi(y)}{|x-y|}\,dy
$$
($\psi$ is a vector field and the equation is interpreted
component-wise). Let $u(x)$ be defined as $u(x)=\curl\psi(x)$, i.e.
\begin{equation}\label{bs}
u(x)=\frac1{4\pi}\int_{\erre^3}\frac{\xi(y)\times(x-y)}{|x-y|^3}\,dy.
\end{equation}
If $\Div\xi=0$, then also $\Div\psi=0$ and $\Div u=0$, and this implies
also $\curl\curl\psi=-\Delta\psi$. Therefore $\curl u=\xi$, i.e.\ $u$ is
the divergence-free velocity field associated to $\xi$. The equality
\eqref{bs} is the \textsl{Biot-Savart law}.

In order to give a probabilistic representation of this formula,
it is necessary to give a representation of the solution of the Poisson
equation and of its derivatives.

\subsection{A probabilistic representation for the Newtonian potential}

In this section we study a probabilistic representation of the solution
of the Poisson equation. The deterministic regularity results are
classical (see for example Gilbarg and Trudinger \cite{GiTr} and Ziemer
\cite{Zie}), so we will focus on the probabilistic formula.

Let $f:\erre^3\to\erre$ be an integrable function. We define the
Newtonian potential with density $f$ as
$$
Nf(x)=\frac1{4\pi}\int_{\erre^3}\frac1{|x-y|}f(y)\,dy.
$$
If $f$ is regular and with compact support, $Nf$ is a solution of the
Poisson equation.

Let $A=\frac12\Delta$, it is well known that $A$ generates, on the space
$C_0(\erre^3)$ of all continuous functions vanishing at infinity, the
strongly continuous semigroup
$$
P_tf(x)=\E[f(x+W_t)] \qquad x\in\erre^3,\ t\ge0,\quad f\in C_0(\erre^3),
$$
where $(W_t)_{t\ge0}$ is a 3D-standard Brownian motion. The resolvent
of $A$ can be written as
$$
\big((A-\lambda I)^{-1}f\big)(x)
=\int_0^\infty\e^{-\lambda t}\E[f(x+W_t)]\,dt,
\qquad f\in C_0(\erre^3),
$$
so that we can argue that the integral
\begin{equation}\label{probpoisson}
\int_0^\infty\E[f(x+W_t)]\,dt.
\end{equation}
converges to $A^{-1}f(x)= 2Nf(x)$ (indeed, at this stage, we do not know
if $A$ is invertible).

As a first step, we find some conditions on $f$ in such a way that
formula \eqref{probpoisson} produces a solution of the Poisson equation.

\begin{proposition}
Let $f\in L^p(\erre^3)\cap L^q(\erre^3)$, with $1\le
p<\frac32<q<\infty$. Then the integral in \eqref{probpoisson} is
convergent for all $x\in\erre^3$ and is equal to $2Nf(x)$. Moreover
$Nf\in C_0(\erre^3)$ and
$$
\|Nf\|_\infty\le C_{p,q}(\|f\|_p+\|f\|_q).
$$
\end{proposition}

\begin{proof}
For every $r>1$, by H\"older inequality,
\begin{equation}\label{tellepi}
\E|f(x+W_t)|
=\frac1{(2\pi t)^{3/2}}\int_{\erre^3}|f(x+y)|\e^{-\frac1{2t}|y|^2}\,dy
\le C_rt^{-\frac3{2r}}\|f\|_r,
\end{equation}
so that, by using the above inequality with $r=p$ and $r=q$ and by
integrating by time,
$$
\int_0^\infty\!\!\!\!\!\E|f(x+W_t)|\,dt
\le\int_0^1\!\!\!\E|f(x+W_t)|\,dt+\int_1^\infty\!\!\!\!\!\E|f(x+W_t)|\,dt
\le C(\|f\|_p+\|f\|_q).
$$
This will prove also the final inequality, once the other properties are
verified. The integral in \eqref{probpoisson} is equal to $2Nf(x)$ since
(we can use Fubini theorem because of the previous inequality)
\begin{eqnarray*}
\int_0^\infty\E[f(x+W_t)]\,dt]
&=&\int_{\erre^3}f(x+y)\int_0^\infty\frac1{(2\pi t)^{3/2}}\e^{-\frac1{2t}|y|^2}\,dt\,dy\\
&=&\int_{\erre^3}\frac1{2\pi|y|}f(x+y)\,dy
 = 2Nf(x).
\end{eqnarray*}
We know from Gilbarg and Trudinger \cite{GiTr}, that $f\in L^q(\erre^3)$
implies, by Sobolev embeddings, $Nf\in C(\erre^3)$. The behaviour at
infinity is less standard, so we give a probabilistic proof of it. Thus
let us show that $Nf\in C_0(\erre^3)$. Indeed, for each $R>0$,
$$
\int_0^\infty\!\!\!\E[f(x+W_t)]\,dt
= \int_0^\infty\!\!\!\E f(x+W_t)I_{\{|W_t|>R\}}\,dt
 +\int_0^\infty\!\!\!\E f(x+W_t)I_{\{|W_t|\le R\}}\,dt
$$
and, in order to show that $Nf(x)$ converges to $0$ as $|x|\to\infty$,
we will prove that the first term converges to $0$, uniformly in $x$, as
$R\to\infty$, and the second term converges to $0$ as $|x|\to\infty$ for
each $R>0$.

For the first term the claim is true since, as in \eqref{tellepi},
$$
\sup_{x\in\erre^3}\E|f(x+W_t)|I_{\{|W_t|>R\}}
\le C(\|f\|_p+\|f\|_q)(t^{-3/2p}I_{[1,\infty)}(t)+t^{-3/2q}I_{[0,1)}(t))
$$
and 
$$
\sup_{x\in\erre^3}\E|f(x+W_t)|I_{\{|W_t|>R\}}
\le Ct^{-3/2}\|f\|_p\bigl(\int_{|y|>R}\e^{-\frac1{2t}|y|^2}\bigr)^{1/p'}
\longrightarrow0
$$
as $R\to\infty$. As regards the second term, we can proceed as in
\eqref{tellepi} and bound the term $\E|f(x+W_t)|I_{\{|W_t|\le R\}}$ with
$$
C\bigl(t^{-\frac3{2p}}\|f(y)I_{\{|y-x|\le R\}}\|_pI_{[1,\infty)}(t)
          +t^{-\frac3{2q}}\|f(y)I_{\{|y-x|\le R\}}\|_qI_{[0,1)}(t)\bigr),
$$
so that, after the integration in time, the above
term converges to $0$, since $f\in L^p(\erre^3)\cap L^q(\erre^3)$.
\end{proof}

In the second step, we study the derivatives of $Nf$. Notice that, for a
regular $f$, Bismut-Elworthy formula (see for example \cite{ElLi}) gives
$$
D_{x_i}\E[f(x+W_t)]
=\frac1t\E[f(x+W_t)W_t^i].
$$
In this simple case, with the Brownian motion, such formula can be
easily checked by means of the Gaussian density.

As in the previous proposition, one could expect that, under
suitable conditions, it is possible to write the derivatives of $Nf$
with the probabilistic representation suggested by the formula
above. Indeed, this is the case, as the following proposition shows.

\begin{proposition}\label{velrep}
Let $f\in L^p(\erre^3)\cap L^q(\erre^3)$ for some $1\le p<\frac32<3<q<+\infty$.
Then $\nabla Nf\in C_0(\erre^3)$ and for each $x\in\erre^3$,
\begin{equation}\label{repD}
2D_{x_i}Nf(x)
=\int_0^\infty\frac1t\E[f(x+W_t)W_t^i],\qquad i=1,\,2,\,3.
\end{equation}
Moreover
\begin{equation}\label{inrepD}
\|\nabla Nf\|_\infty\le C_{p,q}\bigl(\|f\|_p+\|f\|_q\bigr)
\end{equation}
\end{proposition}

\begin{proof}
By H\"older inequality,
\begin{eqnarray}\label{lprepD}
\frac1t\E|f(x+W_t)W^i_t|
& = &\frac{C}{t^{5/2}}\int_{\erre^3}\!\!f(x+y)y^i\e^{-\frac1{2t}|y|^2}\notag\\
&\le&\frac{C}{t^{5/2}}\|f\|_p\sqrt{t}t^{\frac{3}{2p'}}\\
&\le&C\|f\|_pt^{-\frac12-\frac3{2p}}\notag
\end{eqnarray}
and, as in the proof of the previous proposition, the time integral is
finite and bounded with respect to $x$, by the assumptions on $p$ and
$q$. Moreover it can be easily seen, by the same arguments used in the
previous proposition, that formula \eqref{repD} and inequality
\eqref{inrepD} hold and that $\nabla Nf\in C_0(\erre^3)$.
\end{proof}

In the last step, we study the second derivatives of the Newtonian
potential. The regularity of the following theorem is based on the
classical Schauder estimates.

\begin{proposition}\label{velgrad}
Let $f\in L^p(\erre^3)\cap C^\alpha_b(\erre^3)$, with $1\le
p<\frac32$. Then $Nf\in C^{2,\alpha}_b(\erre^3)\cap C_0(\erre^3)$,
$$
\|Nf\|_{C^{2,\alpha}_b(\erre^3)}
\le \widetilde C\bigl(\|f\|_{L^p(\erre^3)}
          +\|f\|_{C^\alpha_b(\erre^3)}\bigr)
$$
and $Nf$ is the unique solution of the Poisson equation in
$C_0(\erre^3)\cap C^2(\erre^3)$.
\end{proposition}

\begin{proof}
From the previous proposition, we know that $Nf\in
C^1_b(\erre^3)$. Bismut-Elworthy formula gives us
\begin{equation}\label{derseconde}
D_{x_ix_j}\E f(x+W_t)
=\frac2t\E[(D_{x_i}\psi)(x+W_{\frac{t}2})W^j_{\frac{t}2}],
\end{equation}
where $\psi(x)=\E f(x+W_{\frac{t}2})$. Hence, in order to show that
$$
D_{x_ix_j}Nf(x)=\int_0^\infty\frac1t\E[(D_{x_i}\psi)(x+W_{\frac{t}2})W^j_{\frac{t}2}]\,dt
$$
holds, it is sufficient to show that \eqref{derseconde} is integrable in
time in the interval $[0,\infty)$.

First, by the Bismut-Elworthy formula, we see that
$$
D_{x_i}\psi(x)=\frac2t\E[f(x+W_{\frac{t}2})W^i_{\frac{t}2}]
$$
and, by \eqref{lprepD}, that
\begin{equation}\label{psiinf}
\|D_{x_i}\psi\|_\infty\le Ct^{-\frac12-\frac3{2p}}\|f\|_p.
\end{equation}
Moreover, since $f\in C^\alpha_b(\erre^3)$,
\begin{eqnarray}\label{psihol}
|D_{x_i}\psi(y)-D_{x_i}\psi(x)|
& = &\frac2t\E|f(y+W_{\frac{t}2})-f(x+W_{\frac{t}2})|\cdot|W^i_{\frac{t}2}|\notag\\
&\le&Ct^{-\frac12}[f]_\alpha|x-y|^\alpha.
\end{eqnarray}
Now we show that \eqref{derseconde} is integrable in time. By
\eqref{psiinf}
$$
\frac2t\E|(D_{x_i}\psi)(x+W_{\frac{t}2})W^j_{\frac{t}2}|
\le Ct^{-\frac32-\frac3{2p}}\|f\|_p\E|W^j_{\frac{t}2}|
\le Ct^{-1-\frac3{2p}}\|f\|_p
$$
and \eqref{derseconde} is integrable in $[1,\infty)$. By \eqref{psihol}
it follows that
\begin{eqnarray*}
\frac2t\E|(D_{x_i}\psi)(x+W_{\frac{t}2})W^j_{\frac{t}2}|
& = & \frac2t\E|\bigl[(D_{x_i}\psi)(x+W_{\frac{t}2})-(D_{x_i}\psi)(x)\bigr]
W^j_{\frac{t}2}|\\
&\le& Ct^{-\frac32}\E|W_{\frac{t}2}|^\alpha|W^i_{\frac{t}2}|\\
&\le& Ct^{-1+\frac\alpha2}[f]_\alpha
\end{eqnarray*}
and \eqref{derseconde} is integrable in $[0,1)$.

In conclusion, the probabilistic representation formula for the second
derivatives holds and
$$
\|D_{x_ix_j}Nf\|_\infty\le C\bigl(\|f\|_p+[f]_\alpha\bigr).
$$
By Schauder's theory, since $Nf\in C^2_b(\erre^3)$ and $f\in
C^\alpha_b(\erre^3)$, it follows that $Nf\in C^{2,\alpha}_b(\erre^3)$
and
$$
\|Nf\|_{C^{2,\alpha}_b(\erre^3)}
\le C\bigl(\|f\|_p+\|f\|_{C^\alpha_b}\bigr)
$$
(see for example Lunardi \cite{Lun}). Moreover, $Nf$ solves the Poisson
equation (Lemma 4.2 of Gilbarg and Trudinger \cite{GiTr}) and the
solution is unique by the maximum principle.
\end{proof}

\subsection{A probabilistic representation for the Biot-Savart law}

We apply now the theory developed in the previous section. The following
theorem, which is actually a mere corollary of the above results, is
nothing but the well known \textsl{Biot-Savart law}.

\begin{theorem}\label{biot}
Let $\xi\in L^p(\erre^3,\erre^3)\cap C^\alpha_b(\erre^3,\erre^3)$, with
$1\le p<\frac32$ and $0<\alpha<1$. There is a unique $u\in
C^{1,\alpha}_b(\erre^3,\erre^3)\cap C_0(\erre^3,\erre^3)$ such that
$$
\curl u=\xi,\qquad \Div u=0
$$
and such solution is given by the following formula
$$
u(x)=\frac12\int_0^\infty\frac1t\E[\xi(x+W_t)\times W_t]\,dt,
\qquad x\in\erre^3,
$$
where $(W_t)_{t\ge0}$ is a standard 3D-Brownian motion.
\end{theorem}

\begin{proof}
The probabilistic formula derives from Proposition \ref{velrep} and the
regularity of $u$ from Propositions \ref{velrep} and \ref{velgrad}. We
prove the uniqueness of the representation: since $\Div u=0$, we have
$u=\curl\psi$, where $\psi$ is the stream function. Now, by the maximum
principle, the unique solution of the problem
$$
\Delta u=0,\qquad\qquad u\to0\quad\text{as }|x|\to\infty
$$
is $u\equiv0$.
\end{proof}

Since we are interested in the time evolution of the vector fields, it
is appropriate to give a time-dependent version of the previous
theorem. We recall that the spaces $\ualfa$ and $\ualfam$  have been
defined in \eqref{ualfa} and \eqref{ualfam}, the spaces $\valfapi$ and
$\valfapil$ have been defined in \eqref{valfapi} and \eqref{valfapil}.

\begin{corollary}\label{biotsavart}
Let $\alpha\in(0,1)$ and $1\le p<\frac32$. The map
$\bs:\valfapi\to\ualfa$, defined as
$$
\bs(\xi)(t,x)=\frac12\int_0^\infty\frac1s\E[\xi(t,x+W_s)\times W_s]\,ds,
$$
is linear bounded and $\|\bs\|\le \widetilde C$, where $\widetilde C$ is
the constant, independent of $T$, appearing in Proposition \ref{velgrad}.

Moreover, if $L,M>0$ are constant such that $M\ge\widetilde CL$, then
the map $\bs:\valfapil\to\ualfam$ is linear bounded.
\end{corollary}


\section{The representation map}\label{lipschitz}

The section is devoted to the study of the properties of the
representation map $\ns$, defined as
$$
\ns(u)(t,x)=\E[U^{x,t}_t\psi(X^{x,t}_t)]+\int_0^t\E[U^{x,t}_tg(t-s,X^{x,t}_s)]\,ds,
$$
where $\psi=\psi(x)$, $g=g(t,x)$ and $X^{x,t}_s$ are the Lagrangian
paths, defined in \eqref{lagpath}, and $U^{x,t}_s$ are the deformation
matrices, defined in \eqref{defmat}.

In the first part, some regularity properties of the Lagrangian paths
and of the deformation matrices are obtained. In the second part we show
that $\ns$ maps the space $\ualfa$ in $\valfapi$ (for the definition of
the spaces, see \eqref{ualfa} and \eqref{valfapi}). Finally, in the
third part, we prove that $\ns$ is Lipschitz-continuous from $\ualfa$ to
$\valfapi$.

\subsection{Regularity of the Lagrangian paths}

In this section we study some regularity properties of the Lagrangian
paths
$$
\left\{\begin{array}{ll}
 dX^x_s=u(s,X^x_s)\,ds+\sqn\,dW_s,&\qquad s\in[0,T],\\
 X_0^x=x,\end{array}\right.
$$
and of the deformation matrices
$$
\left\{\begin{array}{ll}
 dU^x_s=U^x_s\D(s,X^x_s)\,ds,&\qquad s\in[0,T],\\
 U_0^x=I,\end{array}\right.
$$
where $u\in C([0,T];C^1_b(\erre^3,\erre^3))$ and $\D\in
C([0,T];C^\alpha_b(\erre^3,\erre^{3\times 3}))$ are given.
Notice that both equations have unique strong solutions. Hence, fixed a
3D Brownian motion $((W_s)_{s\ge0},(\F_s)_{s\ge0})$ on the
probability space $(\Omega,\F,\Pb)$, for
each $x\in\erre^3$ there is a process $(X^x_s,U^x_s)_{s\ge0}$ that
solves the corresponding equations, and the solution is unique up to
indistinguishability. The equations can be solved path-wise, choosing
the $\omega\in\Omega$ for which $s\to W_s(\omega)$ is a continuous
function. Hence, the statements of this section are true for all such
$\omega$s, independently of $x$ and $s$. First define
$$
\|v\|_{\infty,s}=\sup_{0\le r\le s}\|v(r)\|_\infty.
$$
\begin{lemma}
Assume $u\in C([0,T];C^1_b(\erre^3,\erre^3))$. Then
\begin{equation}\label{laghol}
|X_s^x-X^y_s|\le|x-y|\e^{s\|\nabla u\|_{\infty,s}},
\qquad s\ge0,\ x,y\in\erre^3.
\end{equation}
Moreover, if $\Div u=0$, then for all $s\ge 0$ and $\omega\in\Omega$,
the map 
$$
x\in\erre^3\mapsto X^x_s(\omega)\in\erre^3
$$
is a diffeomorphism, the determinant of its Jacobian is everywhere equal
to $1$ and
\begin{equation}\label{lagleb}
\int_{\erre^3}\varphi(X^x_s(\omega))\,dx=\int_{\erre^3}\varphi(x)\,dx
\qquad\varphi\in L^1(\erre^3).
\end{equation}
\end{lemma}
\begin{proof}
First we prove \eqref{laghol}. By easy computations,
$$
|X_s^x-X_s^y|
\le |x-y|+\|\nabla u\|_{\infty,s}\int_0^s|X_r^x-X_r^y|\,dr
$$
and applying Gronwall's lemma, we can conclude.

Using Theorem 4.6.5 of Kunita \cite{Kun} (actually the assumption of
H\"older continuity on $u$ is useless for our aim, since we deal with an
additive noise, see also Theorem 4.1.1 of Busnello \cite{Bus}) one can
easily deduce that $x\mapsto X^x_t$ is a diffeomorphism and the
determinant of its Jacobian is constant. Moreover, since $\Div u=0$,
the determinant of its Jacobian is equal to $1$ for all times, so that,
by a change of variables and a density argument, also \eqref{lagleb} can
be deduced.
\end{proof}

\begin{lemma}\label{ustima}
Assume $u\in C([0,T];C^1_b(\erre^3,\erre^3))$ and $\D\in
C([0,T];C^\alpha_b(\erre^3,\erre^{3\times 3}))$. Then
$$
|U^x_s|\le\e^{s\|\D\|_{\infty,s}}\qquad x\in\erre^3,\ s\in[0,T]
$$
and for $x,y\in\erre^3$ and $s\in[0,T]$,
$$
|U^x_s-U^y_s|
\le s\e^{2s\|\D\|_{\infty,s}+s\|\nabla
u\|_{\infty,s}}[\D]_{\alpha,s}|x-y|^\alpha.
$$
\end{lemma} 
\begin{proof}
The first property derives from \eqref{Ubound}. As regards the second,
from \eqref{Ubound} and \eqref{laghol} we get
\begin{eqnarray*}
|U^x_s-U^y_s|
&\le& \int_0^s\|\D(r)\|_\infty|U^x_r-U^y_r|\,dr
     +\int_0^s\e^{s\|\D\|_{\infty,s}}|\D(r,X^x_r)-\D(r,X^y_r)|\,dr\\
&\le& \|\D\|_{\infty,s}\int_0^s|U^x_r-U^y_r|\,dr
     +s\e^{s\|\D\|_{\infty,s}+s\|\nabla u\|_{\infty,s}}[\D]_{\alpha,s}|x-y|^\alpha
\end{eqnarray*}
and, by Gronwall's lemma, the claim follows.
\end{proof}

Let ${\mathcal B}_b(\erre^3,\erre^3)$ be the space of all bounded
Borel-measurable functions and define the operator
$$
Q_s\varphi(x)=\E[U^x_s\varphi(X^x_s)],\qquad x\in\erre^3.
$$

\begin{lemma}\label{stimaqt}
Let $s\ge0$, then
\begin{enumerate}
\item[1.] $Q_s\in{\mathcal L}({\mathcal B}_b(\erre^3,\erre^3))$ and
$\|Q_s\|_{{\mathcal L}({\mathcal B}_b)}\le \e^{s\|\D\|_{\infty,s}}$
\item[2.] $Q_s\in{\mathcal L}(C^\alpha_b(\erre^3,\erre^3))$ and
$\|Q_s\|_{{\mathcal L}(C^\alpha_b)}\le \e^{2s\|\D\|_{\infty,s}+s\|\nabla
u\|_{\infty,s}}(1+s[\D]_{\alpha,s})$
\end{enumerate}
Moreover, if $\Div u=0$, then
\begin{enumerate}
\item[3.] $Q_s\in{\mathcal L}(L^p(\erre^3,\erre^3))$ and
$\|Q_s\|_{{\mathcal L}(L^p)}\le \e^{s\|\D\|_{\infty,s}}$
\end{enumerate}
\end{lemma}
\begin{proof}
First property is an obvious consequence of the previous lemma. About
the second, using the two lemmas above,
\begin{eqnarray*}
|\E[U^x_s\varphi(X^x_s)-U^y_s\varphi(X^y_s)]|
&\le& \E|U^x_s-U^y_s|\cdot|\varphi(X^x_s)|
     +\E|U^y_s|\cdot|\varphi(X^x_s)-\varphi(X^y_s)|\\
&\le& (s[\D]_{\alpha,s}+1)\e^{2s\|\D\|_{\infty,s}+s\|\nabla
u\|_{\infty,s}}\|\varphi\|_{C^\alpha_b}|x-y|^\alpha.
\end{eqnarray*}
Finally, assume $\Div u=0$. Using \eqref{lagleb}, H\"older inequality
and the previous lemma, we get
$$
\int_{\erre^3}|Q_s\varphi(x)|^p
\le\e^{ps\|\D\|_{\infty,s}}\E\int_{\erre^3}|\varphi(X^x_s)|^p
\le\e^{ps\|\D\|_{\infty,s}}\|\varphi\|_p^p.
$$
\end{proof}

\subsection{Definition of the representation map}

Here we prove that $\ns$ maps $\ualfa$ in $\valfapi$. Before proving such
claim, we need some preliminary definitions and results. For each
$u\in\ualfa$, consider for all $x\in\erre^3$ and $t\in[0,T]$ the
Lagrangian paths
\begin{equation}\label{lagpath}
\left\{\begin{array}{ll}
 dX^{x,t}_s=-u(t-s,X^{x,t}_s)\,ds+\sqn\,dW_s,&\qquad s\in[0,t],\\
 X_0^{x,t}=x,\end{array}\right.
\end{equation}
and the deformation matrices
\begin{equation}\label{defmat}
\left\{\begin{array}{ll}
 dU^{x,t}_s=U^{x,t}_s\D_u(t-s,X^{x,t}_s)\,ds,&\qquad s\in[0,T],\\
 U_0^{x,t}=I,\end{array}\right.
\end{equation}
where $\D_u=\nabla u$ or $\D_u=\frac12(\nabla u+\nabla u^T)$.

\begin{lemma}
Let $u\in\ualfa$ and $\psi\in C_b(\erre^3,\erre^3)\cap
L^p(\erre^3,\erre^3)$. The function
$$
(s,t)\in\{0\le s\le t\le T\}\mapsto\E[U^{\cdot,t}_s\psi(X^{\cdot,t}_s)]\in
L^p(\erre^3,\erre^3)\cap C_b(\erre^3,\erre^3)
$$
is continuous with respect to both variables.
\end{lemma}
\begin{proof}
First we show the continuity in $C_b$. If $0\le s\le t\le T$ and $0\le
r\le v\le T$, with $t\le v$, then for each $x\in\erre^3$,
\begin{eqnarray}\label{contstima}
\lefteqn{\lvert\E[U^{x,t}_s\psi(X^{x,t}_s)]-\E[U^{x,v}_r\psi(X^{x,v}_r)]\rvert\le}\notag\\
\quad
&\le& \E|U^{x,t}_s-U^{x,v}_s|\cdot|\psi(X^{x,t}_s)|
      +\E|U^{x,v}_s|\cdot|\psi(X^{x,t}_s)-\psi(X^{x,v}_s)|\\
&   & +\E|U^{x,v}_s-U^{x,v}_r|\cdot|\psi(X^{x,v}_s)|
      +\E|U^{x,v}_r|\cdot|\psi(X^{x,v}_s)-\psi(X^{x,v}_r)|.\notag
\end{eqnarray}
In order to estimate the different terms of the above inequality, we see
that from the equations
$$
|U^{x,v}_s-U^{x,v}_r|
=   |\int_s^rU^{x,v}_\s\D_u(v-\s,X^{x,v}_\s)\,d\s|
\le \e^{v\|\D_u\|_{\infty,v}}\|\D_u\|_{\infty,v}|r-s|,
$$
and
\begin{equation}\label{fintest}
|X^{x,v}_r-X^{x,v}_s|
\le \|u\|_\infty|s-r|+\sqn|W_r-W_s|.
\end{equation}
Moreover
\begin{eqnarray*}
|X^{x,t}_s-X^{x,v}_s|
&\le&\int_0^s|u(t-\s,X^{x,t}_\s)-u(v-\s,X^{x,v}_\s)|\,d\s\\
&\le&\int_0^s\|u(t-\s)-u(v-\s)\|_\infty
     +\int_0^s\|\nabla u(v-\s)\|_\infty|X^{x,t}_\s-X^{x,v}_\s|
\end{eqnarray*}
and, by Gronwall's lemma,
$$
|X^{x,t}_s-X^{x,v}_s|
\le\e^{v\|\nabla u\|_{\infty,v}}\int_0^s\|u(t-\s)-u(v-\s)\|_\infty\,d\s.
$$
Finally,
\begin{eqnarray*}
|U^{x,t}_s-U^{x,v}_s|
&\le& \int_0^s|U^{x,t}_\s|\cdot|\D_u(t-\s,X^{x,t}_\s)-\D_u(t-\s,X^{x,v}_\s)|\,d\s\\
&   &+\int_0^s|U^{x,t}_\s|\cdot|\D_u(t-\s,X^{x,v}_\s)-\D_u(v-\s,X^{x,v}_\s)|\,d\s\\
&   &+\int_0^s|\D_u(v-\s,X^{x,v}_\s)|\cdot|U^{x,t}_\s-U^{x,v}_\s|\,d\s\\
&\le&\e^{v\|\nabla u\|_{\infty,v}}\|\D_u\|_{\infty,v}
        \sup_{[0,v]}|X^{x,t}_\s-X^{x,v}_\s|\\
&   &+v\e^{v\|\nabla u\|_{\infty,v}}\sup_{[0,v]}\|\D_u(t-\s)-\D_u(v-\s)\|_\infty\\
&   &+\|\nabla u\|_{\infty,v}\int_0^s|U^{x,t}_\s-U^{x,v}_\s|\,d\s\\
&\le&A(t,v)+C\int_0^s|U^{x,t}_\s-U^{x,v}_\s|\,d\s,
\end{eqnarray*}
where $\E A(t,v)\to0$ as $|t-v|\to0$, and by Gronwall's lemma,
$$
|U^{x,t}_s-U^{x,v}_s|\le A(t,v)\e^{Cv}.
$$
Using the above estimates in \eqref{contstima}, it is easy to show
continuity with values in $C_b$. In order to show continuity in $L^p$,
we remark that the above estimates ensure convergence for all
$x\in\erre^3$, so that we need only to show uniform integrability. To
this end, notice that, by a change of variables,
\begin{eqnarray*}
\int_{|x|\ge K}|\E[U^{x,t}_s\psi(X^{x,t}_s)]|^p
&\le& C\E\int_{X^{\cdot,t}_s(B_K^c)}|\psi(y)|^p\,dy\\
&\le& C\int_{|y|\ge\frac{K}2}|\psi|^p\Pb[|X^{x,t}_s-x|\le\frac{K}2]\\
&   &+C\|\psi\|_p^p\Pb[|X^{x,t}_s-x|\ge\frac{K}2],\\
&\le& C\int_{|y|\ge\frac{K}2}|\psi|^p
     +C\|\psi\|_p^p\Pb[|X^{x,t}_s-x|\ge\frac{K}2],
\end{eqnarray*}
where $C=T\e^{T\|\nabla u\|_{\infty,T}}$, and, because of
\eqref{fintest}, for $K\to\infty$, the above quantity converges to $0$
independently of $s$, $t$.
\end{proof}

Now it is possible to prove the above mentioned result on the map $\ns$.

\begin{proposition}\label{nsmap}
Given $1\le p<\frac32$ and $0<\alpha<1$, let $\psi\in
C^\alpha_b(\erre^3,\erre^3)\cap L^p(\erre^3,\erre^3)$ and
$g\in\valfapi$, then $\ns$ maps $\ualfa$ in $\valfapi$ and
\begin{equation}\label{NSbound}
\|\ns(u)(t)\|_{C^\alpha_b\cap L^p}
\le\e^{3t\|\nabla u\|_{\infty,t}}(1+t\|\nabla u\|_{C^\alpha_b})
       \bigl(\|\psi\|_{C^\alpha_b\cap L^p}+\int^t_0\|g(s)\|_{C^\alpha_b\cap L^p}\,ds\bigr)
\end{equation}
\end{proposition}
\begin{proof}
First, $\ns(u)\in C^\alpha_b\cap L^p$ follows by Lemma \ref{stimaqt},
moreover also estimate \eqref{NSbound} can be easily deduced. Finally,
from the previous lemma it follows that
$$
t\mapsto\ns(u)(t)\in C^\alpha_b\cap L^p
$$
is continuous. 
\end{proof}

\subsection{Lipschitz continuity of the representation map}

Let $g\in\valfapi$ and $\psi\in C^\alpha_b(\erre^3,\erre^3)\cap
L^p(\erre^3,\erre^3)$, and consider the map
$$
\ns:\ualfa\longrightarrow\valfapi
$$
defined in the previous section. The aim of the present section is to
show that such map is locally Lipschitz-continuous. In order to do this,
we will use Girsanov formula. First we rewrite $\ns$ in a more
appropriate form, namely
$$
\ns(u)(t,x)=\E[F_{t,u}(X^{x,t,u})],
$$
where for each trajectory $w\in C([0,T];\erre^3)$,
$$
F_{t,u}(w)=V^{t,u}_t(w)\psi(w_t)+\int_0^tV^{t,u}_s(w)g(t-s,w_s)\,ds
$$
and $V^{t,u}(w)$ is the solution of the following differential equation
$$
\left\{\begin{array}{ll}\dot V^{t,u}_s=V^{t,u}_s\D_u(t-s,w_s),
    &\qquad s\le t,\\
  V^{t,u}_0(w)=I.\\\end{array}\right.
$$
Notice that $U^{x,t,u}_s(\omega)=V^{t,u}_s(X^{x,t,u}(\omega))$, for each
$\omega\in\Omega$, and we have made an explicit reference to the
dependence from $u$ in the Lagrangian paths $X^{x,t,u}$ and in the
deformation matrices $U^{x,t,u}$.

By Girsanov formula, we have
$$
\E[F_{t,u}(X^{x,t,u})]=\E[Z_t^{x,t,u}F_{t,u}(x+\sqn W_\cdot)],
$$
where
$$
Z_s^{x,t,u}
=\exp\bigl[\frac1{\sqn}\int_0^s\langle u(t-r,x+\sqn W_r),dW_r\rangle
     -\frac1{4\nu}\int_0^s|u(t-r,x+\sqn W_r)|^2\,dr\bigr],
$$
with $s\le t$, so that for each $u$,
\begin{eqnarray*}
\ns(u)(t,x)
&=&\E[Z^{x,t,u}_tV^{t,u}_t(x+\sqn W_\cdot)\psi(x+\sqn W_t)]+\\
& &+\int_0^t\E[Z^{x,t,u}_tV^{t,u}_s(x+\sqn W_\cdot)g(t-s,x+\sqn W_t)]\,ds.
\end{eqnarray*}
Using this representation, we will prove the following proposition.

\begin{proposition}\label{nslip}
Given $1\le p<\frac32$ and $0<\alpha<1$, let $\psi\in
C^\alpha_b(\erre^3,\erre^3)\cap L^p(\erre^3,\erre^3)$ and
$g\in\valfapi$ and set
$$
\ep_0=\|\psi\|_{C^\alpha_b\cap L^p}+\int_0^T\|g(s)\|_{C^\alpha_b\cap L^p}\,ds.
$$
For each $u$, $v\in\ualfam$,
$$
\sup_{t\le T}\|\ns(u)-\ns(v)\|_{C_b^\alpha\cap L^p}
\le C(\nu,p)C_M(T)\ep_0\sup_{t\le T}\|u-v\|_{C^{1,\alpha}_b},
$$
where $C(\nu,p)$ is a constant depending only on $p$ and $\nu$ and
$\lim_{T\to0}C_M(T)=0$.
\end{proposition}

The proof of the above proposition will be carried on using the
subsequent lemmas. In order to make the explanations easier, we introduce
the following notations. We define $\Delta_{xy}f=f(x)-f(y)$ for any
function $f$. Notice that
\begin{equation}\label{delta}
\Delta_{xy}(fg)=(\Delta_{xy}f)g(x)+f(y)(\Delta_{xy}g).
\end{equation}
If the functions depends on two variables, we define $\Delta_{uvxy}$ as
$\Delta_{uv}\Delta_{xy}$ and, by applying twice the above formula,
\begin{eqnarray}\label{ddelta}
\Delta_{uvxy}(fg)
&=& \Delta_{uv}[(\Delta_{xy}f)g(\cdot,x)+f(\cdot,y)(\Delta_{xy}g)]\notag\\
&=& (\Delta_{uvxy}f)g(u,x)+[\Delta_{xy}f(v)][\Delta_{uv}g(x)]\\
& &+[\Delta_{uv}f(y)][\Delta_{xy}g(u)]+f(v,y)(\Delta_{uvxy}g).\notag
\end{eqnarray}
\begin{lemma}
Let $u$, $v\in\ualfam$, then for each $w$, $w'\in C([0,T];\erre^3)$ and
for all $s\le t$,
\begin{eqnarray*}
&&|V^{t,u}_s(w)|\le\e^{tM},\\
&&|\Delta_{uv}V^{t,\cdot}_s(w)|\le t\e^{2tM}\|u-v\|_{C^1_b},\\
&&|\Delta_{ww'}V^{t,u}_s(\cdot)|\le 2Mt\e^{2tM}\|w-w'\|_\infty^\alpha,\\
&&|\Delta_{uvww'}V^{t,\cdot}_s(\cdot)|
   \le(1+3tM)t\e^{3tM}\|w-w'\|_\infty^\alpha\|u-v\|_{C^{1,\alpha}_b}.
\end{eqnarray*}
\end{lemma}
\begin{proof}
The proofs of these properties are similar, we just give the
proof of the last one. Indeed, using formula \eqref{ddelta},
\begin{eqnarray*}
\frac{d}{ds}(\Delta_{uvww'}V^{t,\cdot}_s(\cdot))
&=&\Delta_{uvww'}\bigl(\frac{d}{ds}V^{t,\cdot}_s(\cdot)\bigr)
 = \Delta_{uvww'}(V^{t,\cdot}_s(\cdot)\D_\cdot(t-s,\cdot))\\
&=& [\Delta_{uvww'}V^{t,\cdot}_s(\cdot)]\D_u(t-s,w_s)
   +V^{t,v}_s(w')[\Delta_{ww'}\D_{u-v}(t-s,\cdot)]\\
& &+[\Delta_{ww'}V^{t,v}_s(\cdot)]\D_{u-v}(t-s,w_s)
   +[\Delta_{uv}V^{t,\cdot}_s(w')][\Delta_{ww'}\D_u(t-s,\cdot)]
\end{eqnarray*}
so that, by using the other inequalities of this lemma,
\begin{eqnarray*}
|\Delta_{uvww'}V^{t,\cdot}_s(\cdot)|
&\le& M\int_0^s|\Delta_{uvww'}V^{t,\cdot}_r(\cdot)|\,dr
     +\|w-w'\|_\infty^\alpha\|u-v\|_{C^{1,\alpha}_b}\int_0^s|V^{t,v}_r(w')|\,dr\\
&   &+\|u-v\|_{C^{1,\alpha}_b}\int_0^s|\Delta_{ww'}V^{t,v}_r(\cdot)|\,dr
     +M\|w-w'\|_\infty^\alpha\int_0^s|\Delta_{uv}V^{t,\cdot}_r(w')|\,dr\\
&\le& M\int_0^s|\Delta_{uvww'}V^{t,\cdot}_r(\cdot)|\,dr
     +(1+3tM)s\e^{2tM}\|w-w'\|_\infty^\alpha\|u-v\|_{C^{1,\alpha}_b}
\end{eqnarray*}
and, by the Gronwall's lemma, the inequality follows.
\end{proof}
Using the previous lemma and formulas \eqref{delta} and \eqref{ddelta},
we can easily deduce similar properties for the functional $F$.
\begin{lemma}\label{sullaf}
Let $u$, $v\in\ualfam$, then for each $w$, $w'\in C([0,T];\erre^3)$, and
for all $t\in[0,T]$,
\begin{eqnarray*}
&&|F_{t,u}(w)|\le\e^{tM}[|\psi(w_t)|+\int_0^t|g(t-s,w_s)|\,ds]\\
&&|\Delta_{uv}F_{t,\cdot}(w)|\le t\e^{2tM}\|u-v\|_{C^1_b}[|\psi(w_t)|+\int_0^t|g(t-s,w_s)|\,ds]\\
&&|\Delta_{ww'}F_{t,u}(\cdot)|\le(1+2tM)\e^{2tM}\ep_0\|w-w'\|_\infty^\alpha\\
&&|\Delta_{uvww'}F_{t,\cdot}(\cdot)|\le (2+3tM)t\e^{3tM}\ep_0\|w-w'\|_\infty\|u-v\|_{C^{1,\alpha}_b},
\end{eqnarray*}
where $\ep_0=\|\psi\|_{C^\alpha_b}+\int_0^t\|g(s)\|_{C^\alpha_b}\,ds$.
\end{lemma}
Finally, we estimate the same quantities on the process $Z$.
\begin{lemma}\label{sullaz}
Let $u$, $v\in\ualfam$ and $q\ge2$. Then for all $s\le t$,
\begin{eqnarray*}
&&\E|Z^{x,t,u}_s|^q\le C\e^{Ct^{q/2}M^q},\\
&&\E|\Delta_{uv}Z^{x,t,\cdot}_s|^q\le Ct^{q/2}\e^{CM^qt^{q/2}}\|u-v\|_{C_b}^q,\\
&&\E|\Delta_{xy}Z^{\cdot,t,u}_s|^q\le Ct^{q/2}M^q\e^{CM^qt^{q/2}}|x-y|^{\alpha q},\\
&&\E|\Delta_{uvxy}Z^{\cdot,t,\cdot}_s|^q
    \le Ct^{3q/2}M^{2q}\e^{CM^qt^{q/2}}|x-y|^{\alpha q}\|u-v\|_{C^{1,\alpha}_b}^q,
\end{eqnarray*}
where $C=C(q,\nu)$ is a constant depending only on $q$ and $\nu$.
\end{lemma}
\begin{proof}
From the definition, we see that $Z^{x,t,u}_s$ solves
$$
\left\{\begin{array}{ll}
dZ^{x,t,u}_s=\frac1{\sqn}Z^{x,t,u}_su(t-s,x+\sqn W_s)\,dW_s,
&\qquad s\le t,\\
Z^{x,t,u}_0=1.
\end{array}\right.
$$
Again, the proofs of the four inequalities are similar, we prove only
the last one. By applying formula \eqref{ddelta}, we get
\begin{eqnarray*}
\lefteqn{d_s(\Delta_{uvxy}Z^{\cdot,t,\cdot}_s)
=\Delta_{uvxy}(d_sZ^{\cdot,t,\cdot}_s)=}\\
&=& \frac1\sqn\Bigl[(\Delta_{uvxy}Z^{\cdot,t,\cdot}_s)u(t-s,Y^x_s)\,dW_s
   +Z^{y,t,v}_s\Delta_{xy}[u(t-s,Y^\cdot_s)-v(t-s,Y^\cdot_s)]\,dW_s\\
& &+(\Delta_{xy}Z^{\cdot,t,v}_s)[u(t-s,Y^x_s)-v(t-s,Y^x_s)]\,dW_s
   +(\Delta_{uv}Z^{y,t,\cdot}_s)[\Delta_{xy}u(t-s,Y^\cdot_s)]\,dW_s\Bigr],
\end{eqnarray*}
where, for the sake of briefness, we have set $Y^x_s=x+\sqn W_s$. By the
Burkholder, Davis and Gundy inequality,
\begin{eqnarray*}
\E|\Delta_{uvxy}Z^{\cdot,t,\cdot}_s|^q
&\le& C\Bigl[M^q\E[\int_0^s|\Delta_{uvxy}Z^{\cdot,t,\cdot}_r|^2dr]^{\frac{q}2}
     +\|u-v\|_{C^\alpha_b}^q|x-y|^{\alpha q}\E[\int_0^s|Z^{y,t,v}_r|^2dr]^{\frac{q}2}\\
&   &+\|u-v\|_{C_b}^q\E[\int_0^s|\Delta_{xy}Z^{\cdot,t,v}_r|^2dr]^{\frac{q}2}
     +M^q|x-y|^{\alpha q}\E[\int_0^s|\Delta_{uv}Z^{y,t\cdot}_r|^2dr]^{\frac{q}2}\Bigr],
\end{eqnarray*}
so that, by using the H\"older inequality and the other inequalities of
this lemma, we get
$$
\E|\Delta_{uvxy}Z^{\cdot,t,\cdot}_s|^q
\le CM^qs^{\frac{q}2-1}\int_0^s\E|\Delta_{uvxy}Z^{\cdot,t,\cdot}_r|^q\,dr
   +CM^{2q}t^qs^{q/2}\e^{CM^qt^{q/2}}|x-y|^{\alpha q}\|u-v\|_{C^{1,\alpha}_b}^q.
$$
Finally, using the Gronwall's lemma, we obtain the required inequality.
\end{proof}
We are now able to prove the main result of this section.
\begin{proof}[Proof of Proposition \ref{nslip}]
Let $u$, $v\in\ualfam$. We start with the estimates in $C_b$ and $L^p$.
Using formula \eqref{delta} and H\"older inequality, we get, for each
$x\in\erre^3$ and $t>0$,
\begin{eqnarray}\label{cibiellepi}
|[\Delta_{uv}\ns(\cdot)](t,x)|
& = & |\E[\Delta_{uv}\bigl(Z^{x,t,\cdot}_tF_{t,\cdot}(Y^x)\bigr)]|\notag\\
&\le& C(q)\Bigl[\bigl(\E|\Delta_{uv}Z^{x,t,\cdot}_t|^{q'}\bigr)^{1/{q'}}\bigl(\E|F_{t,u}(Y^x)|^q\bigr)^{1/q}\\
&   &+\bigl(\E|Z^{x,t,v}_t|^{q'}\bigr)^{1/{q'}}\bigl(\E|\Delta_{uv}F_{t,\cdot}(Y^x)|^q\bigr)^{1/q}\Bigr]\notag
\end{eqnarray}
where $q\ge1$, $q'$ is the H\"older conjugate exponent of $q$ and
we have set $Y^x_s=x+\sqn W_s$. Using the estimates in Lemma
\ref{sullaf} and in Lemma \ref{sullaz}, and the inequality above with
$q=2$, we obtain the estimate in the $C_b$ norm,
$$
\sup_{t\le T}\|\Delta_{uv}\ns(\cdot)\|_{C_b}
\le C\ep_0(T+\sqrt{T})\e^{(CM^2+2M)T}\|u-v\|_{C^{1,\alpha}_b}.
$$
Using  again Lemma \ref{sullaf} and \ref{sullaz} and the inequality
\eqref{cibiellepi} above, with $q=p$, we can obtain the estimate in the 
$L^p$ norm,
$$
\sup_{t\le T}\|\Delta_{uv}\ns(\cdot)\|_{L^p}^p
\le C\ep_0^p(T^p+T^{p/2})\e^{2TM+CM^{p'}t^{p'/2}}\|u-v\|_{C^{1,\alpha}_b}.
$$

To conclude the proof, we need the estimate in the $C^\alpha_b$ norm. For
all $x$, $y\in\erre^3$ and $t>0$, by applying formula \eqref{ddelta} we get
\begin{eqnarray*}
|\Delta_{uvxy}\ns(\cdot)(t,\cdot)|
&\le& \E[\Delta_{uvxy}(Z^{\cdot,t,\cdot}F_{t,\cdot}(Y^\cdot))]\\
&\le& \E\Bigl[(\Delta_{uvxy}Z^{\cdot,t,\cdot}_t)F_{t,u}(Y^x)
     +Z^{y,t,v}_t[\Delta_{uvxy}F_{t,\cdot}(Y^\cdot)]\\
&   &+(\Delta_{uv}Z^{y,t,\cdot}_t)[\Delta_{xy}F_{t,u}(Y^\cdot)]
     +(\Delta_{xy}Z^{\cdot,t,v}_t)[\Delta_{uv}F_{t,\cdot}(Y^x)]\Bigr]
\end{eqnarray*}
Using the inequalities in Lemma \ref{sullaf} and Lemma \ref{sullaz}, it
follows that
$$
|\Delta_{uvxy}\ns(\cdot)(t,\cdot)|
\le C\ep_0(\sqrt{T}+T+MT^{3/2}+MT^2)\e^{3TM+CTM^2}|x-y|^\alpha\|u-v\|_{C^{1,\alpha}_b}.
$$
\end{proof}


\bibliographystyle{amsplain}

\begin{thebibliography}{9}
%
\bibitem{AlBe} \textsc{S.~Albeverio, Ya.~I.~Belopolskaya}, Work in preparation.
%
\bibitem{ArKh} \textsc{V.~I.~Arnold, B.~A.~Khesin}, \textsl{Topological
Methods in Hydrodynamics}, Appl. Math. Sci. \textbf{125}, Springer
Verlag, Berlin, 1998.
%
\bibitem{BKM} \textsc{J.~T.~Beale, T.~Kato, A.~Majda}, \textsl{Remarks on
the breakdown of smooth solutions for the 3D Euler equations},
Comm. Math. Phys. \textbf{94} (1984), 61--66.
%
\bibitem{Ben} \textsc{M.~Ben-Artzi}, \textsl{Global Solutions of
Two-Dimensional Navier-Stokes and Euler Equations},
Arch. Anal. Rat. Mech. \textbf{128} (1994), 329--358.
%
\bibitem{Bus} \textsc{B.~Busnello}, \textsl{A probabilistic approach to
the two-dimensional Navier-Stokes equations}, Ann. Probab. \textbf{27}
no. 4 (1999), 1750--1780.
%
\bibitem{Can} \textsc{M.~Cannone}, \textsl{Viscous flows in Besov
spaces}, Advances in mathematical fluid mechanics (Paseky, 1999), 1--34,
Springer, Berlin, 2000.
%
\bibitem{Cho} \textsc{A.~J.~Chorin}, \textsl{Vorticity and turbulence},
Springer-Verlag, New York, 1994.
%
\bibitem{Con} \textsc{P.~Constantin}, \textsl{Geometric statistics in
turbulence}, SIAM Review \textbf{36}(I) (1994), 73--98.
%
\bibitem{ElLi} \textsc{K.~D.~Elworthy, X.~M.~Li}, \textsl{Formulae for
the derivatives of heat semigroups}, J. Funct. Anal. \textbf{125}
(1994), 252--286.
%
\bibitem{EMPS} \textsc{R.~Esposito, R.~Marra, M.~Pulvirenti,
C.~Sciarretta}, \textsl{A stochastic Lagrangian picture for the three
dimensional Navier-Stokes equations}, Comm. Partial
Diff. Eq. \textbf{13}(12) (1988), 1601--1610.
%
\bibitem{EsPu} \textsc{R.~Esposito, M.~Pulvirenti},
\textsl{Three-dimensional stochastic vortex flows}, Math. Methods
Appl. Sci. II (1989), 431--445.
%
\bibitem{Fre} \textsc{M.~I.~Freidlin}, \textsl{Functional Integration
and Partial Differential Equations}, Princeton Univ. Press, Princeton,
1985.
%
\bibitem{Gie} \textsc{J.~S.~Giet}, Ph.D. Thesis, Institute Eli\'e
Cartan, Universit\'e H. Poincar\'e (Nancy I).
%
\bibitem{GiTr}\textsc{D.~Gilbarg, N.~S.~Trudinger}, \textsl{Elliptic
Partial Differential Equations of Second Order}. Second edition.
Grundlehren der Mathematischen Wissenschaften, \textbf{224}.
Springer-Verlag, Berlin, 1983.
%
\bibitem{Gli} \textsc{Y.~Gliklikh}, \textsl{Global analysis in
mathematical physics. Geometric and stochastic methods},
Translated from the 1989 Russian original and with Appendix F by Viktor
L. Ginzburg. Applied Mathematical Sciences,
\textbf{122}. Springer-Verlag, New York, 1997.
%
\bibitem{Kah} \textsc{C.~S.~Kahane}, \textsl{The Feynman-Kac formula for
a system of parabolic equations}, Czech. Math. J. \textbf{44}(4) (1994),
579--602.
%
\bibitem{Kry} \textsc{N.~V.~Krylov}, \textsl{On Kolmogorov's equations
for finite dimensional diffusions}, in \textsl{Stochastic PDE's and
Kolmogorov Equations in Infinite Dimensions}, G.~Da~Prato Ed., Lecture
Notes in Mathematics \textbf{1715}, Springer, 1999.
%
\bibitem{Kun} \textsc{H.~Kunita}, \textsl{Stochastic flows and
stochastic differential equations}. Cambridge Studies in Advanced
Mathematics, \textbf{24}. Cambridge University Press, Cambridge, 1990.
%
\bibitem{LeSz} \textsc{Y.~Le~Jan, A.~S.~Sznitman}, \textsl{Stochastic
cascades and $3$-dimensional Navier-Stokes equations},
Prob. Th. Rel. Fields \textbf{109}(3) (1997), 343--366.
%
\bibitem{Lun} \textsc{A.~Lunardi}, \textsl{Analytic semigroups and
optimal regularity in parabolic problems}. Progress in Nonlinear
Differential Equations and their Applications, \textbf{16}.
Birkh\"auser Verlag, Basel, 1995.
%
\bibitem{Rap} \textsc{D.~L.~Rapoport}, \textsl{Closed form integration
of the Navier Stokes equations with stochastic differential geometry},
Hadronic J. \textbf{22} (1999), no. 5, 577--605.
%
\bibitem{vWa} \textsc{W.~von Wahl}, \textsl{The equations of
Navier-Stokes and abstract parabolic equations}, Aspects of Mathematics,
E8. Friedr. Vieweg \& Sohn, Braunschweig,  1985.
%
\bibitem{Zie} \textsc{W.~P.~Ziemer}, \textsl{Weakly differentiable
functions}. Graduate Texts in Mathematics, \textbf{120}.
Springer Verlag, New York, 1989.
%
\end{thebibliography}

\end{document}